\documentclass[12pt]{amsart}

\usepackage{amscd}
\usepackage{amsfonts}
\usepackage{amssymb}
\usepackage{amsmath}
\usepackage{latexsym}

\newcommand{\la}{\lambda}
\newcommand{\al}{\alpha}
\newcommand{\be}{\beta}
\newcommand{\ga}{\gamma}

\newcommand{\f}{\varphi}
\newcommand{\h}{\phi}

\newcommand{\abs}[1]{\vert #1\vert}
\newcommand{\ov}{\overline}
\newcommand{\Ll}{\mathcal{L}}

\newcommand{\GL}{\mathop{\rm GL}\nolimits}

\newcommand{\ra}{{\rightarrow}}
\newcommand{\lra}{{\longrightarrow}}
\newcommand{\Null}{\mathop{\rm Null}\nolimits}
\newcommand{\Aut}{\mathop{\rm Aut}\nolimits}

\newcommand{\dev}{\mathop{\rm dev}\nolimits}

\newcommand{\CC}{\mathbb{C}} 
 
\newcommand{\RR}{\mathbb{R}}
\newcommand{\ZZ}{\mathbb{Z}}

\numberwithin{equation}{section}

\newtheorem*{ta}{Theorem A}
\newtheorem*{tb}{Theorem B}

\newtheorem*{tc}{Theorem C}
\newtheorem*{ca'}{Corollary A$_1$}
\newtheorem{pr}{Proposition}[section]
\newtheorem{te}{Theorem}[section]
\newtheorem{co}{Corollary}[section]
\newtheorem{lm}{Lemma}[section]

\theoremstyle{definition}
\newtheorem{de}{Definition}[section] 
\newtheorem{re}{Remark}[section]

\begin{document}

\title[Geometric flow on compact l.c.K. manifolds]
{Geometric flow on 
compact locally conformally K\"ahler manifolds}
\author{Yoshinobu Kamishima}
\address{Department of Mathematics, Tokyo Metropolitan 
University,\newline
Minami-Ohsawa 1-1, Hachioji, Tokyo 192-0397, Japan}
\email{kami@comp.metro-u.ac.jp}

\author{Liviu Ornea}
\address{University of Bucharest, Faculty of Mathematics\newline
14 Academiei str.,
70109 Bucharest, Romania}
\email{lornea@imar.ro}
\date{\today}
\thanks{The second author is a member of EDGE, Research Training Network
HPRN-CT-2000-00101, supported by The European Human Potential Programme}

\begin{abstract}
We study two kinds of transformation groups of a compact locally conformally
K\"ahler (l.c.K.) manifold. 
First we study compact l.c.K. manifolds
by means of the existence of holomorphic l.c.K. flow 
($i.e.,$ a conformal, holomorphic flow with respect 
to the Hermitian metric.) We characterize the structure of the
compact l.c.K. manifolds with parallel Lee form.
Next, we introduce the Lee-Cauchy-Riemann (LCR) transformations as a class of
diffeomorphisms preserving the specific $G$-structure of l.c.K.
manifolds. We show that compact l.c.K. manifolds with
parallel Lee form admitting a $\CC^*$ flow of LCR transformations are rigid: it is holomorphically isometric to
a Hopf manifold with parallel Lee form. 

\end{abstract}

\maketitle


\section{Introduction}

Let $(M,g,J)$ be a connected, complex Hermitian manifold of complex 
dimension $n\geq 2$. We denote its  
fundamental $2$-form by $\omega$; it is  defined by $\omega(X,Y) =g(X,JY)$. If there exists a 
real $1$-form $\theta$ satisfying the integrability condition 

\begin{equation*}
d\omega=\theta\wedge\omega\quad \text{with}\; d\theta=0
\end{equation*}
then $g$ is said to be  a \emph{locally conformally K\"ahler} (l.c.K.) metric.
A complex manifold $M$ endowed with a l.c.K. metric is called a l.c.K. manifold. The conformal class of 
a l.c.K. metric $g$ is said to be a l.c.K. structure on $M$. 
The closed $1$-form $\theta$ is called \emph{the Lee form} and it encodes
the geometric properties of such a manifold. The vector field $\theta^\sharp$,
defined by $\theta(X)=g(X,\theta^\sharp)$, 
is called the Lee field.

The purpose of this paper is to study two kinds of transformation groups of a l.c.K. manifold $(M,g,J)$. 
We first consider $\mathop{\Aut}_{l.c.K.}(M)$, the  group of all conformal,
holomorphic diffeomorphisms. We discuss its properties in \S 2.
A holomorphic vector field $Z$ on $(M,g,J)$ generates a $1$-dimensional
complex Lie group $\mathcal C$. (The universal covering group of 
$\mathcal C$ is $\CC$.) We call $\mathcal C$ a holomorphic flow
on $M$.

\begin{de}\label{parallellck}
If a holomorphic flow $\mathcal C$ (resp. holomorphic vector field $Z$)
belongs to $\mathop{\Aut}_{l.c.K.}(M)$
 (resp. Lie algebra of $\mathop{\Aut}_{l.c.K.}(M)$), then 
$\mathcal C$ (resp. $Z$) is said to be a \emph{holomorphic l.c.K. flow}
(resp. \emph{holomorphic l.c.K. vector field}).
\end{de}

A nontrivial subclass of l.c.K. manifolds is formed by those $(M,g,J)$ having
parallel Lee form w.r.t. the
Levi-Civita connection $\nabla^g$ ($i.e.$ $\nabla^g\theta=0$).
We observe that a compact non-K\"ahler l.c.K. manifold $(M,g,J)$
with parallel Lee form $\theta$ supports a
holomorphic vector field $Z=\theta^\sharp-iJ\theta^\sharp$ 
which generates holomorphic isometries of
$g$. (Compare \cite{Va},\cite{Va1},\cite{DO}.)
We shall prove that the converse is also true:

\begin{ta}
Let $(M,g,J)$ be a compact, connected, l.c.K. non-K\"ahler
manifold, of complex dimension at least $2$. If $\mathop{\Aut}_{l.c.K.}(M)$ contains a holomorphic l.c.K. flow, then there exists a metric with 
parallel Lee form in the
conformal class of $g$.
\end{ta}

\begin{ca'}
With the same hypothesis, $M$  admits a l.c.K. metric with parallel Lee form
if and only if it admits a  holomorphic l.c.K. flow. 
\end{ca'}

In \S 3,
we discuss the existence of l.c.K. metrics with parallel Lee form
on the Hopf manifold. (Compare with \cite{GO}). 
Let $\Lambda=(\la_1, \ldots ,\la_n)$ with the $\la_i$'s complex numbers 
satisfying  $0<\abs{\la_n}\leq \cdots \leq \abs{\la_1}<1$.
By a \emph{primary Hopf manifold} $M_\Lambda$ of type $\Lambda$
we mean the compact quotient manifold 
of $\CC^n-\{0\}$ by a subgroup $\Gamma_\Lambda$ generated by the
transformation $(z_1,\ldots ,z_n)\mapsto (\la_1z_1,\ldots ,\la_nz_n)$.
Note that a primary Hopf manifold of type $\Lambda$ of complex dimension $2$ is  a primary Hopf surface of K\"ahler rank $1$. 
We prove the following: 

\begin{tb} 
The primary Hopf manifold $M_\Lambda$ of type $\Lambda$
supports a l.c.K. metric with parallel Lee form.
\end{tb} 

More generally, we prove the existence
of a l.c.K. metric with parallel Lee form on the Hopf manifold
(cf. Theorem \ref{secondary}).

In the second half of the paper we adopt the viewpoint of $G$-structure
theory in order to 
study a non-compact, non-holomorphic, transformation group of a 
compact l.c.K. manifold $(M,g,J)$. 
Locally, the $2$-form $\omega$ defines the 
real $1$-forms $\theta$, $\theta\circ J$ and $(n-1)$ complex  $1$-forms $\theta^\al$ and their 
conjugates $\bar\theta^\al$, 
where $\theta\circ J$ is called the \emph{anti-Lee form} and is defined by 
$\theta\circ J(X)=\theta(JX)$. 
We consider the group $\mathop{\Aut}_{LCR}(M)$ of transformations of $M$
preserving the structure of unitary coframe fields 
$\mathcal{F}=\{\theta,\theta\circ J, \theta^1, \ldots ,\theta^{n-1}, \bar\theta^1, \ldots, \bar\theta^{n-1}\}$. More precisely,
an element $f$ of $\mathop{\Aut}_{LCR}(M)$  is called a 
\emph{Lee-Cauchy-Riemann} (LCR) transformation if  
it satisfies the equations:
\begin{equation*}
\begin{split}
f^*\theta&=\theta,\\
 f^*(\theta\circ J)&=\lambda\cdot(\theta\circ 
J),\\
f^*\theta^{\al}&=\sqrt \lambda\cdot\theta^\be U^{\al}_{\be}+
(\theta\circ J)\cdot v^{\al},\\
f^*{\bar \theta}^{\al}&=\sqrt \lambda\cdot
{\bar\theta}^\be\ov{U}^\al_\be+(\theta\circ J)\cdot\ov{v}^{\al}.
\end{split}
\end{equation*}
Here $\lambda$ is a positive, smooth function, and
$v^\al\in \CC$, $U^\al_\be\in{\rm U}(n-1)$ are smooth functions.
Obviously, if
$\mathrm{I}(M,g,J)$ is the group of holomorphic isometries,
then both $\mathop{\Aut}_{l.c.K.}(M)$ and $\mathop{\Aut}_{LCR}(M)$
contain $\mathrm{I}(M,g,J)$.

As the main result of this part we exhibit the rigidity of compact l.c.K. manifolds under the existence of a 
non-compact LCR flow:

\begin{tc}
Let $(M,g,J)$ be a compact, connected, l.c.K. non-K\"ahler  manifold of
complex dimension at least $2$, with  parallel Lee form $\theta$. 
Suppose that $M$ admits a closed subgroup $\CC^*=S^1\times \RR^+$
of Lee-Cauchy-Riemann transformations whose $S^1$
subgroup induces the Lee field $\theta^\sharp$.
Then $M$ is holomorphically isometric, up to scalar multiple of the metric,
to the primary Hopf  manifold $M_\Lambda$ of type $\Lambda$.
\end{tc}

\section{Locally conformally K\"ahler transformations}

\begin{pr}\label{compact}
Let $(M,g,J)$ be a compact l.c.K. manifold with $\dim_\CC M\geq 2$. Then $\mathop{\Aut}_{l.c.K.}(M)$ is a compact Lie group.
\end{pr}

\begin{proof}
Note that $\mathop{\Aut}_{l.c.K.}(M)$ is a closed Lie subgroup
in the group of all conformal diffeomorphisms of $(M,g)$. 
If $\mathop{\Aut}_{l.c.K.}(M)$
were noncompact, then
by the celebrated result of Obata and Lelong-Ferrand (\cite{Ob}, \cite{Le}), $(M,g)$ would be 
conformally equivalent
with the sphere $S^{2n}$, $n\geq 2$. Hence $M$ would be  simply connected.
It is well known that a compact simply connected
l.c.K. manifold is conformal to a K\"ahler manifold
(cf. \cite{DO}), which is impossible because the sphere $S^{2n}$ has no K\"ahler structure. 

\end{proof}

From now on, we shall suppose that the l.c.K. manifolds we work with 
are compact, non-K\"ahler and, moreover, the Lee form is not identically zero 
at any point of the manifold. In particular, these manifolds are not simply 
connected (cf. \cite{DO}).
Given a l.c.K. manifold $(M,g,J)$, let $\tilde M$ be the universal covering space of $M$, let $p:\tilde M\rightarrow M$ be the 
canonical projection and denote also by $J$ the lifted complex structure on $\tilde M$. We 
can associate to the fundamental $2$-form
$\omega$ a canonical K\"ahler form on $\tilde M$ as follows. Since the lee
form $\theta$ 
is closed, its lift to $\tilde M$ is exact, hence $p^*\theta=d\tau$ for some 
smooth function $\tau$ 
on $\tilde M$. We put $h=e^{-\tau}\cdot p^*g \;(\text{resp.}\;
 \Omega=e^{-\tau}\cdot p^*\omega).$
It is easy to check that $d\Omega=0$, thus $h$ is a K\"ahler metric on
 $(\tilde M,J)$. In particular \emph{ $g$ is locally conformal to the K\"ahler
 metric $h$} (compare with \cite{DO} and the bibliography therein).
Let $f\in \mathop{\Aut}_{l.c.K.}(M)$. By definition,
$f^*\omega=e^{\lambda}\cdot \omega$ for some function 
$\lambda$ on $M$. Differentiate this equality
to yield that
$(f^*\theta-\theta+d\lambda)\wedge \omega=0$. As $\omega$ is
nondegenerate and ${\rm dim}_{\CC}\ M>1$, $f^*\theta=\theta+d\lambda$. 
Since $p^*\theta=d\tau$, for any lift  $\tilde f$ of $f$ to $\tilde M$ 
we have $d\tilde f^*\tau=d(\tau +p^*\la)$, 
thus $-\tilde f^*\tau + p^*\lambda=-\tau+c$ for some constant $c$.
We can write $\tilde f^*\Omega=e^c\cdot\Omega.$
If $c\neq 0$, $\tilde f$ is a holomorphic homothety w.r.t. $h$;
when  $c=0$,  $\tilde f$ will be an isometry.

We denote by $\mathcal H(\tilde M,\Omega,J)$ the  group of all holomorphic,
 homothetic transformations of the universal 
cover $\tilde M$ w.r.t. the K\"ahler structure $(h,J)$. 
If $f_1,\, f_2\in \mathcal H(\tilde M,\Omega,J)$, there exists
some constant $\rho(f_i)$ ($i=1,2$) satisfying 
 $f_i^*\Omega=\rho(f_i)\cdot \Omega$ as above. It is easy to check
 that $\rho(f_1\circ f_2)=\rho({f_1})\cdot \rho(f_2)$.
We obtain a continuous homomorphism:

\begin{equation}\label{assign}
\rho:\mathcal H(\tilde M,\Omega,J)\lra \RR^+.
\end{equation}
Let $\pi_1(M)$ be the fundamental group of $M$.
Then we note that $\pi_1(M)\subset \mathcal H(\tilde M,\Omega,J)$.
For this, if $\gamma\in\pi_1(M)$, then
$\gamma^*\Omega=e^{-\gamma^*\tau}\cdot \gamma^*p^*\omega
=e^{-\gamma^*\tau}\cdot p^*\omega
=e^{-\gamma^*\tau+\tau}\cdot\Omega$.
Since $\Omega$ is a K\"ahler form $(n\geq 2)$,
$e^{-\gamma^*\tau+\tau}$ must be constant $\rho(\gamma)$.

Let $\mathcal C$ be a holomorphic l.c.K. flow on $M$.
If we denote $\tilde{\mathcal C}$ a lift of $\mathcal C$ to $\tilde M$,
then $\tilde{\mathcal C}\subset \mathcal H(\tilde M,\Omega,J)$.
If $V$ is a vector field which generates a one-parameter subgroup
of $\tilde{\mathcal C}$, then so does $JV$ such as
$V$ and $JV$ together generate $\tilde{\mathcal C}$.
We define a smooth function $s:\tilde M\rightarrow \RR$ to be
$s(x)=\Omega(JV_x,V_x)$. Since
$\tilde{\mathcal C}$ centralizes each element $\ga$
of $\pi_1(M)$, it follows that
$s(\gamma x)=\Omega(JV_{\gamma x},V_{\gamma x})=
\Omega(\gamma_*JV_x,\gamma_*V_x)=\rho(\gamma)s(x).$ If every element
$\gamma$ satisfies that $\rho(\gamma)=1$, $i.e.,$ $\gamma^*\Omega=\Omega$,
then $\pi_1(M)$ acts as holomorphic isometries of $h$ so that
$\Omega$ would induce a K\"ahler structure on $M$. By our hypothesis,
this does not occur. There exists at least one  element $\gamma$
such that $\rho(\gamma)\neq 1$. In particular, we note that:

\begin{equation}\label{nonconst}
\mbox{The function}\ s\ \mbox{is not constant on}\  \tilde M.\ \ \ \ \ \ \ \
\ \ \
\end{equation}
On the other hand, we prove the following lemma. (The proof of the lemma
is almost same as that of \cite{kamicomp}.)

\begin{lm} \label{ess}
$\rho(\tilde{\mathcal C})=\RR^+$, $i.e.,$
the group $\tilde{\mathcal C}$
acts by holomorphic, non-trivial
ho\-mo\-the\-ties w.r.t. the K\"ahler metric
$h$ on $\tilde M$.
\end{lm}

\begin{proof}
Suppose that
$\rho(\tilde{\mathcal C})=\{1\}$. Then $\tilde{\mathcal C}$ leaves
$\Omega$ invariant. As $\{V,JV\}$ generates $\tilde{\mathcal C}$,
it follows that $\Ll_V\Omega=\Ll_{JV}\Omega=0$.
In particular, $Vs=(JV)s=0$.
For any distribution $D$ on $\tilde M$, denote by $D^\perp$ the
orthogonal complement to $D$ w.r.t. the metric $h$ where
$h(\tilde X,\tilde Y)=\Omega(J\tilde X,\tilde Y)$.
Since $0=(\Ll_V\Omega)(JV,\tilde X)=
V\Omega(JV,\tilde X)-\Omega([V,JV],\tilde X)-
\Omega(JV,[V,\tilde X])$, if $\tilde X\in \{V,JV\}^{\perp}$, then
$\Omega(JV,[V,\tilde X])=0$, similarly $\Omega(V,[JV,\tilde X])=0$.
The equality
\begin{equation*}
\begin{split}
0=3d\Omega(\tilde X, V,JV)&=\tilde X\Omega(V,JV)-
V\Omega(\tilde X,JV)+JV\Omega(\tilde X,V)\\
&\ \ -\Omega([\tilde X,V],JV)
-\Omega([V,JV],\tilde X)-\Omega([JV,\tilde X],V)
\end{split}
\end{equation*} implies that $\tilde X\Omega(V,JV)=0$, $i.e.,$
$\tilde Xs=0$ for any $\tilde X\in \{V,JV\}^{\perp}$.
Therefore, $s$ becomes constant, being a contradiction
to \eqref{nonconst}.

\end{proof}

\subsection{The submanifold $W$ and its pseudo-Hermitian structure}
\label{sub-ps}
As ${\rm Ker}\ \rho$ has one dimension, denote by
$-J\xi$ the vector field  whose one-parameter subgroup
$\{\psi_t\}_{t\in \RR}$ acts as holomorphic isometries on $\tilde M$.

\begin{equation}\label{rho2}
\psi_t^*\Omega =\Omega, \quad t\in \RR.
\end{equation}
Since $-J\xi$ and $\xi$ together generate the group $\tilde{\mathcal C}$,
the $1$-parameter subgroup $\{\f_t\}_{t\in\RR}$
generated by $\xi$ acts as nontrivial
holomorphic homotheties w.r.t. $\Omega$ by Lemma \ref{ess}.
In particular, the group
$\{\f_t\}_{t\in\RR}$ is isomorphic to $\RR$.
Since  $\f_t^*\Omega =\rho(\f_t)\cdot \Omega$
 \ $(t\in \RR, \, \rho(\f_t)\in \RR^+)$ from \eqref{assign} and
$\rho$ is a continuous homomorphism, $\rho(\f_t)=e^{at}$
for some constant $a\neq 0$. We may normalize $a=1$
so that:
\begin{equation}\label{rho1}
\f_t^*\Omega =e^t\cdot \Omega, \quad t\in \RR.
\end{equation}

\begin{lm}\label{proper}
The group $\{\f_t\}_{t\in\RR}$ acts properly and hence freely on
$\tilde M$. In particular, $\xi\neq 0$ everywhere on $\tilde M$.
\end{lm}

\begin{proof}
Recall that $\mathcal C$ lies in $\mathop{\Aut}_{l.c.K.}(M)$ by definition.
As $\mathop{\Aut}_{l.c.K.}(M)$ is a compact Lie group,
its closure $\overline{\mathcal C}$ in $\mathop{\Aut}_{l.c.K.}(M)$ is
also compact and so isomorphic to a $k$-torus $(k\geq 2)$.
Therefore, the lift $H$ of $\overline{\mathcal C}$ to $\tilde M$
acts properly on $\tilde M$.
The lift $H$ is isomorphic to $\RR^\ell\times T^m$ where
$\ell+m=k$. Note that $\ell\geq 1$ because
$\rho$ maps any compact subgroup of $H$
to $\{1\}$, but the group $\{\f_t\}_{t\in\RR}\subset H$ satisfies
$\rho(\{\f_t\})=\RR^+$.  Hence the group
$\{\f_t\}_{t\in\RR}$ has a nontrivial summand in
$\RR^\ell$ which implies that $\{\f_t\}_{t\in\RR}$ is closed in
$H$. Thus, the group $\{\f_t\}_{t\in\RR}$
acts properly on $\tilde M$. If we note that
$\{\f_t\}_{t\in\RR}$ is isomorphic to $\RR$,
then it acts freely on $\tilde M$.

\end{proof}

\begin{pr} \label{s1}
Let $s:\tilde M\rightarrow \RR$ be the smooth map defined as
$s(x)=\Omega(J\xi_x,\xi_x)$. Then $1$ is a regular value of $s$, hence
$s^{-1}(1)$ is a codimension one, regular submanifold of $\tilde M$.
\end{pr}

\begin{proof}
As $\f_t$ is holomorphic,
$s(\f_tx)=\Omega(J\xi_{\f_tx},\xi_{\f_tx})=
\Omega(\f_{t*}J\xi_x,\f_{t*}\xi_x)=e^t\cdot s(x)$.
Hence,
$$\Ll_\xi s=
\lim_{t\rightarrow 0}\frac{\f_t^*s-s}{t}=s.
$$
We note also that
\begin{equation}\label{lieomega}
\Ll_\xi\Omega=\Omega.
\end{equation}
By Lemma \ref{proper},
notice that $\xi\neq 0$ everywhere on $\tilde M$.
Since $s(x)\neq 0$, $s^{-1}(1)\neq \emptyset$.
For $x\in s^{-1}(1)$,
$ds(\xi_x)=(\Ll_\xi s)(x)=s(x)=1.$
This proves that $ds:T_x\tilde M\rightarrow \RR$ is onto and so
$s^{-1}(1)$ is a codimension one smooth regular submanifold of $\tilde M$.

\end{proof}

Let now $W=s^{-1}(1)$. We can prove:

\begin{lm}\label{neighbor}
The submanifold $W$ is connected and
the  map $H:\RR\times W\rightarrow \tilde M$,
defined by $H(t, w)=\f_tw$ is an equivariant diffeomorphism.
\end{lm}

\begin{proof}
Let $W_0$ be a component of $s^{-1}(1)$ and $\RR\cdot W_0$ be the
set
$\{\f_tw\; ;\; w\in W_0, t\in \RR\}$. As $\RR=\{\f_t\}$
acts freely and
$s(\f_tx)=e^ts(x)$, we have $\f_tW_0\cap W_0=\emptyset$ for $t\neq
0$.
Thus $\RR\cdot W_0$ is an open subset of $\tilde M$. We prove that
it
is also closed. Let $\ov{\RR\cdot W_0}$ be the closure of $\RR\cdot
W_0$ in $\tilde M$. We choose a limit point $p=\lim \f_{t_i}w_i\in
\ov{\RR\cdot W_0}$. Then $s(p)=\lim s(\f_{t_i}w_i)=\lim
e^{t_i}s(w_i)=\lim e^{t_i}$. Put $t=\log s(p)$, then $t=\lim t_i$,
so $\f_t^{-1}(p)=\lim \f_{t_i}^{-1}(\lim \f_{t_i}w_i)=\lim w_i$.
Since $s^{-1}(1)$ is regular (\emph{i.e.} closed w.r.t. the
relative topology induced from $\tilde M$), its component $W_0$
is also closed. Hence $\f_t^{-1}p\in W_0$.
Therefore $p=\f_t(\f_t^{-1}p)\in \RR\cdot
W_0$, proving that $\RR\cdot W_0$ is closed in $\tilde M$. In
conclusion, $\RR\cdot W_0=\tilde M$. Now, if $W_1$ is another
component of $s^{-1}(1)$, the same argument shows
 $\RR\cdot W_1=\tilde M$. As $\RR\cdot W_0=\RR\cdot W_1$ and
$s(W_1)=1$, this implies $W_0=W_1$, in other words $W$ is
connected.

\end{proof}

Let $i:W\rightarrow \tilde M$ be the inclusion
and $\pi:\tilde
M\rightarrow W$ be the canonical projection.
Define a $1$-form $\eta$ on $W$ to be
\begin{equation}\label{et}
\eta=i^*\iota_\xi\Omega.
\end{equation}
Here $\iota_\xi$ denotes the interior product with $\xi$.
We have from $\S$ \ref{sub-ps} that:
\begin{equation}\label{diff-vec}
\frac{d\psi_t}{dt}(x)|_{t=0}=-J\xi_x.
\end{equation}Using \eqref{rho2},
$s(\psi_tw)=s(w)=1$ $(w\in W)$ so that the group
$\{\psi_t\}_{t\in\RR}$ leaves $W$ invariant. Hence, the vector field
$-J\xi$ restricts to a vector field $A$ to $W$.
If $\{\psi'_t\}_{t\in\RR}$
is the one-parameter subgroup
generated by $A$, then
\begin{equation}\label{comm-inc}
\psi_t=i\circ \psi'_t.
\end{equation}

\begin{lm} \label{re}
The $1$-form $\eta$ is a contact form on $W$ for which $A$ is the
characteristic vector field $($Reeb field$)$.
\end{lm}

\begin{proof} First note that
$\eta(A_w)=\iota_\xi\Omega(-J\xi_w)=\Omega(J\xi_w,\xi_w)=s(w)=1\ \ (w\in W).$
Moreover, from \eqref{lieomega},
$d\eta=i^*d\iota_\xi\Omega=i^*(d\iota_\xi\Omega+\iota_\xi d\Omega)=
i^*\Ll_\xi\Omega=i^*\Omega$. Hence,
$\eta\wedge d\eta^{n-1}\neq 0$ on $W$ showing that $\eta$ is a
contact form.
Noting \eqref{rho2}, \eqref{comm-inc} and that
both $\f_t$ and $\psi_\theta$
commutes each other, it is easy to see that
\begin{equation}\label{psi-inv}
\begin{split}
{\psi'}_t^*\iota_\xi\Omega&=\iota_\xi\Omega\ \ \ \mbox{on}\ \tilde M.\\
{\psi'}_t^*\eta&=\eta\ \ \ \mbox{on}\ W.
\end{split}
\end{equation} Let $\mathop{\Null}\, \eta=\{X\in TW\ |\ \eta(X)=0\}$ be the
contact subbundle.
Since $\Ll_A\eta(X)=A\eta(X)-\eta([A,X])$ and $\Ll_A\eta=0$
from \eqref{psi-inv}, if $X\in\mathop{\Null}\, \eta$, then
$\eta([A,X])=0$. Moreover,\\
$\displaystyle d\eta(A,X)=\frac 12(A\eta(X)-X\eta(A)-
\eta([A,X]))=0$, which implies that $d\eta(A,X)=0$ for all $X\in TW$,
showing that $A$ is the characteristic vector field.

\end{proof}

Recall that $\RR\rightarrow \tilde M\stackrel{\pi}{\rightarrow}W$ is a
principal fiber bundle with $T\RR=<\xi>$.
By Lemma \ref{neighbor}, each point $x\in \tilde M$
can be described uniquely as $x=\f_tw$.
Using \eqref{comm-inc},
\begin{equation}\label{comm-pro}
\begin{split}
&\pi\circ\psi_\theta(x)=\pi\circ\psi_\theta(\f_tw)=
\pi\circ\f_t(\psi_\theta w)\\
&=\pi\circ i{\psi'}_\theta(w)={\psi'}_\theta(w)={\psi'}_\theta\circ\pi(x),
\end{split}
\end{equation}hence, $\pi_*(-J\xi)=A$.
As $i_*\pi_*X_{x}-X_x=a\cdot \xi_x$ for some function $a$,
using \eqref{et}, $\pi$ maps $\{\xi,J\xi\}^\perp$
isomorphically onto $\mathop{\Null}\, \eta$.
Since $\{\xi,J\xi\}^\perp$ is $J$-invariant,
there exists an almost complex structure $J$ on
$\mathop{\Null}\, \eta$ such that the following diagram is commutative:

\begin{equation}\label{J-com}
\begin{CD}
\{\xi,J\xi\}^\perp@>\pi_*>>{\Null}\, \eta\\
@VVJ V                 @VVJ V\\
\{\xi,J\xi\}^\perp@>\pi_*>>{\Null}\, \eta.
\end{CD}
\end{equation}
\begin{pr}\label{psh}
The pair $(\eta, J)$ is a strictly pseudoconvex, pseudo-Her\-mi\-tian
structure on $\tilde W$.
\end{pr}
\begin{proof} Let $\Psi: {\Null}\, \eta\times {\Null}\, \eta
\ra \RR$ be the bilinear form defined by
$\Psi(X,Y)=d\eta(JX,Y)$. There exist $\tilde X,\ \tilde Y\in
\{\xi,J\xi\}^\perp$ such that
$\pi_*\tilde X=X$, $\pi_*\tilde Y=Y$. Then it is easy to see that
$i_*JX\equiv J\tilde X,\ \ i_*Y\equiv \tilde Y \ \mbox{mod}\ \xi.$
Using $d\eta=i^*\Omega$ as above,
$\Psi(X,Y)=i^*\Omega(JX,Y)=\Omega(J\tilde X,\tilde Y)
=h(\tilde X,\tilde Y),$
hence $\Psi$ is positive definite.
By definition, $\eta$ is strictly pseudoconvex.
Let $\{\xi,J\xi\}^\perp\otimes \CC=B^{1,0}\oplus B^{0,1}$ be
the canonical splitting  of $J$.
Then we prove that
$[B^{1,0}, B^{1,0}]\subset B^{1,0}$.
Let $\tilde X,\tilde Y \in B^{1,0}$.
Since $T^{1,0}\tilde M=\{\xi-{\mathrm i}J\xi\}\oplus B^{1,0}$
and $J$ is integrable on $\tilde M$,
$[\tilde X,\tilde Y]\in T^{1,0}\tilde M$.
Put $[\tilde X,\tilde Y]=a(\xi-{\mathrm i}J\xi)+\tilde Z$
for some  function $a$ and $\tilde Z\in B^{1,0}$.
As $\pi_*(-J\xi)=A$ from \eqref{comm-pro},
$\pi_*([\tilde X,\tilde Y])=a{\mathrm i}A+\pi_*\tilde Z$.
By definition,
$2d\eta(\pi_*\tilde X,\pi_*\tilde Y)=
-\eta([\pi_*\tilde X,\pi_*\tilde Y])=-a{\mathrm i}$.
On the other hand, since $\Omega$ is $J$-invariant,
$\Omega(\tilde X,\tilde Y)=0$ for $\forall\ \tilde X,\tilde Y \in B^{1,0}$.
As above, $i_*\pi_*\tilde X\equiv \tilde X$ $\mbox{mod}\ \xi$, similarly for
$\tilde Y$, we obtain that $d\eta(\pi_*\tilde X,\pi_*\tilde Y)=
\Omega(i_*\pi_*\tilde X, i_*\pi_*\tilde Y)=\Omega(\tilde X,\tilde Y)=0.$
Hence, $a=0$ and so $[\tilde X,\tilde Y]=\tilde Z\in B^{1,0}$.
If we note that $\pi_*:\{\xi,J\xi\}^\perp\otimes 
\CC\ra {\Null}\, \eta\otimes\CC$ is $J$-isomorphic by \eqref{J-com},
then ${\Null}\, \eta\otimes\CC=\pi_*B^{1,0}\oplus \pi_*B^{0,1}$
is the splitting for $J$, in which we have shown
$[\pi_*B^{1,0}, \pi_*B^{1,0}]\subset\pi_* B^{1,0}$.
Therefore $J$ is a complex structure on 
${\Null}\, \eta$.

\end{proof}

Consider the group of pseudo-Hermitian transformations on $(W,\eta,J)$:
\begin{equation}\label{PSH}
\mathrm{PSH}(W,\eta,J)=\{f\in \mathrm {Diff}(W) \ |\
 f^*\eta=\eta, f_*\circ J=J\circ f_*\ \mbox{on}\, {\Null}\ \eta\}.
\end{equation}

\begin{co}\label{psh1}
The characteristic vector field $A$
generates
the subgroup
$\{{\psi'}_t\}_{t\in \RR}$
consisting of pseudo-Her\-mi\-tian transformations.  
\end{co}

\begin{proof}
By \eqref{rho2} and \eqref{psi-inv},
${\psi}_t$ (resp. ${\psi'}_t$)
preserves $\{\xi,J\xi\}^\perp$ (resp. $\mathop{\Null}\, \eta$).
Then 
the equality $\pi\circ\psi_\theta={\psi'}_\theta\circ\pi$
from \eqref {comm-pro} with diagram \eqref{J-com}
implies that
${\psi'_t}_*J=J{\psi'_t}_*$ on
${\Null}\, \eta$.
Therefore
\begin{equation}\label{crtrans}
\{\psi'_t\}_{t\in\RR}\subset {\rm PSH}(W,\eta,J).
\end{equation}

\end{proof}

\subsection*{Proof of Theorem A}\hfill
\subsection{Parallel Lee form} 
Let $Y_{\f_tw}\in T_{\f_tw}\tilde M$ be any vector field.
As $\pi_*Y_{\f_tw}\in T_wW$, 
$i_*\pi_*Y_{\f_tw}-{\f_{-t}}_*Y_{\f_tw}=\lambda \xi_w$
for some function $\lambda$.
Then,
\begin{equation*}
\begin{split}
&\iota_\xi\Omega(i_*\pi_*Y_{\f_tw})=
\Omega(\xi_w,i_*\pi_*Y_{\f_tw})=\Omega(\xi_w,{\f_{-t}}_*Y_{\f_tw})
+\Omega(\xi_w,\lambda \xi_w)\\
&=\f_{-t}^*\Omega({\f_t}_*\xi_w,Y_{\f_tw})=
e^{-t}\Omega(\xi_{\f_tw},Y_{\f_tw})=
e^{-t}\iota_\xi\Omega(Y_{\f_tw}).
\end{split}
\end{equation*}By definition \eqref{et},
\begin{equation}\label{form-eta}
\pi^*\eta=\pi^*i^*\iota_\xi\Omega
=e^{-t}\iota_\xi\Omega, \ \ \mbox{equivalently,}\ \
e^t\pi^*\eta=\iota_\xi\Omega.
\end{equation}
As $\Omega=\Ll_\xi\Omega=d\iota_\xi\Omega$ from \eqref{lieomega},
we obtain that

\begin{equation}\label{eta}
d(e^{t}\pi^*\eta)=\Omega \;\, \text{on}\  \tilde M.
\end{equation}

For the given l.c.K. metric $g$, the K\"ahler metric $h$ is obtained as $h=e^{-\tau}\cdot p^*g$ where
$d\tau=\tilde \theta$. 
As $\omega$ is the fundamental $2$-form of $g$, note
that $\Omega=e^{-\tau}\cdot p^*\omega$.

We now consider on $\tilde M$ the $2$-form:
\begin{equation}\label{bar}
\bar\Theta=2e^{-t}\cdot d(e^{t}\pi^*\eta)\ (=2e^{-t}\cdot\Omega).
\end{equation}
Then $\bar g(X,Y)=\bar\Theta(JX,Y)$ is a l.c.K. metric.
Put $\bar \theta=-dt$. Then, as\\
$d\bar\Theta=-2e^{-t}dt\wedge d(e^{t}\pi^*\eta)=-dt\wedge\bar\Theta$,
so $\bar \theta$ is the Lee form of $\bar g$.

\begin{lm}\label{parallelform}
$\bar \theta$ is parallel w.r.t. $\bar g$
$(\nabla^{\bar g}\bar\theta=0)$.
\end{lm}

\begin{proof}
 First we determine the Lee field $\bar\theta^{\sharp}$.
$(\bar\theta(X)=\bar g(X,\bar\theta^{\sharp}).)$
We start from:
\begin{equation*}
\begin{split}
\bar g(\xi,Y)&=\bar\Theta(J\xi,Y)=2e^{-t}(e^{t}dt
\wedge\pi^*\eta+e^{t}d\pi^*\eta)(J\xi,Y)\\
&=2(dt\wedge \pi^*\eta+d\pi^*\eta)(J\xi,Y)=2(dt\wedge \pi^*\eta)(J\xi,Y)
\end{split}
\end{equation*}
because $A=-\pi_*J\xi$ is the characteristic vector field of
the contact form $\eta$. 
As before, a point $x\in \tilde M$ can be described uniquely
as $\f_tw$ for some $w\in W$. In particular, using
Lemma \ref{neighbor}, the $t$-coordinate of $x$ is $t$.
Noting that
$\psi_\theta(x)=\f_t{\psi_\theta}w$
and ${\psi_\theta} w\in W$, by uniqueness
the $t$-coordinate of $\psi_\theta(x)$, 
$t(\psi_\theta(x))=t$.
From \eqref{diff-vec},
\begin{equation}\label{diff}
dt(-J\xi_x)=dt(\frac{d\psi_\theta}{d\theta}(x)|_{\theta=0})=
\frac {dt}{d\theta}|_{\theta=0}=0.
\end{equation}
The above formula becomes: 
\begin{equation}\label{sharp}
\bar g(\xi,Y)=2(dt\wedge \pi^*\eta)(J\xi,Y)=-dt(Y)\eta(-A)
=dt(Y)=
-\bar\theta(Y)=-\bar g(Y,\bar\theta^\sharp)
\end{equation}
proving that $\bar\theta^\sharp=-\xi$.
Next we observe that the flow $\{\f_s\}_{s\in\RR}$ acts by isometries
w.r.t. $\bar g$. As $\f_s$ is holomorphic, it is enough to prove that each $\f_s$ leaves $\bar \Theta$ invariant. But
\begin{equation*}
\f_s^*\bar \Theta=2e^{-\f_s^*t}d(e^{\f_s^*t}\f_s^*\pi^*\eta)=2e^{-(s+t)}
d(e^{s+t}\pi^*\eta)
=2e^{-t}d(e^t\pi^*\eta)=\bar\Theta.
\end{equation*}
Thus $\mathcal{L}_{\theta^{\sharp}}\bar g=-\mathcal{L}_\xi\bar g=0$. 
Now  we put $\sigma=\bar \theta$ in the equality
$\displaystyle (\mathcal{L}_{\sigma^\sharp}\bar g)(X,Y)+
2d\sigma(X,Y)=2\bar g(\nabla^{\bar g}_X\sigma^\sharp, Y)$,
valid for any $1$-form $\sigma$,  
take into account $d\bar\theta=0$ and obtain $\displaystyle \nabla^{\bar g}
\bar \theta^\sharp=0$ which is equivalent with $\displaystyle \nabla^{\bar g}
\bar \theta=0$, so $\bar\theta$ is parallel w.r.t. $\bar g$ as announced. 

\end{proof}

By  equation \eqref{bar}, $\bar g$ is conformal to the lifted metric $p^*g$:
\begin{equation}\label{co}
\bar\Theta=\mu\cdot p^*\omega\ \ \text{(equivalently}\
\bar g=\mu\cdot p^*g)
\end{equation}
where  $\mu=2e^{-(t+\tau)}:
\tilde M\ra \RR^+$ is a smooth map. We finally prove:

\begin{lm}\label{inv}
$\pi_1(M)$ acts by holomorphic isometries of $\bar g$.
In particular,
$\pi_1(M)$ leaves $\bar\theta$ invariant.
\end{lm}
\begin{proof}
We prove the following two facts:
\begin{itemize}
\item[{\bf 1.}] $\gamma^*\pi^*\eta=\pi^*\eta$ for every $\gamma\in\pi_1(M)$.
\item[{\bf 2.}]  $\gamma^*e^t=\rho(\gamma)\cdot e^t$ where $\rho:\pi_1(M)\ra \RR^+$
is the homomorphism as before.
\end{itemize}
First note that as $\RR=\{\f_t\}$ centralizes $\pi_1(M)$,
$\gamma_*\xi=\xi$ for $\gamma\in\pi_1(M)$.
As $\gamma$ is holomorphic, $\gamma_*J\xi=J\xi$. 
Since $\pi_1(M)$ acts on $\tilde M$ as  holomorphic homothetic transformations,
$(i.e.,\ \gamma^*\Omega=\rho(\gamma)\cdot\Omega)$,
$\pi_1(M)$ preserves $\{\xi,J\xi\}^\perp$.
If we recall that
$\pi_*:\{\xi,J\xi\}^\perp\rightarrow\mathop{\Null}\, \eta$
is isomorphic, then for $X\in\{\xi,J\xi\}^\perp$, $\gamma^*\pi^*\eta(X)
=\eta(\pi_*\gamma_*X)=0$. As $-\pi_*J\xi=A$ is  characteristic,
it follows $\gamma^*\pi^*\eta(J\xi)
=\eta(\pi_*\gamma_*J\xi)=\eta(\pi_*J\xi)=-1$.
This shows that
$\gamma^*\pi^*\eta=\pi^*\eta$ on $\tilde M$. On the other hand,
if we note $\gamma_*\xi=\xi$,
then
\begin{eqnarray*}
\gamma^*(\iota_\xi\Omega)(X)&=&\Omega(\xi,\gamma_*X)=
\Omega(\gamma_*\xi,\gamma_*X)
=\gamma^*\Omega(\xi,X)\\
&=&\rho(\gamma)\cdot\Omega(\xi,X)
=\rho(\gamma)\cdot \iota_\xi\Omega(X)
\end{eqnarray*}
where $\rho(\gamma)$ is a positive constant number.
Applying $\gamma^*$ to $\pi^*\eta=e^{-t}\cdot\iota_\xi\Omega$ from
\eqref{form-eta},
we obtain $\gamma^*e^{-t}\cdot\rho(\gamma)=e^{-t}$.
Equivalently, $\gamma^*e^{t}=\rho(\gamma)\cdot e^{t}$.
This shows {\bf 1} and {\bf 2}.\\
From \eqref{bar},
\begin{eqnarray*}
\gamma^*\bar\Theta&=&\gamma^*(2e^{-t}\cdot d(e^{t}\pi^*\eta))
 =2\rho(\gamma)^{-1}\cdot e^{-t}d(\rho(\gamma)\cdot
e^{t}\gamma^*\pi^*\eta)\\
&=&2e^{-t}\cdot d(e^{t}\pi^*\eta)=\bar\Theta.
\end{eqnarray*}
Since $\bar g(X,Y)=\bar\Theta(JX,Y)$,
$\pi_1(M)$ acts through holomorphic isometries of $\bar g$.
We have that $\bar\theta(Y)=\bar g(Y,\bar\theta^\sharp)
=-\bar g(Y,\xi)$ $(Y\in T\tilde M)$ from \eqref{sharp}.
Then,
\begin{equation*}
\gamma^*\bar\theta(Y)=-\bar g(\gamma_*Y,\xi)
=-\bar g(\gamma_*Y,\gamma_*\xi)=-\bar g(Y,\xi)=\bar\theta(Y).
\end{equation*}

\end{proof}

From this lemma,  the covering map $p:\tilde M\ra M$ induces a
l.c.K. metric $\hat g$ with parallel Lee form $\hat\theta$ 
on $M$ such that $p^*\hat g=\bar g$ and $p^*\hat\theta=\bar\theta$
with $\nabla^{\hat g}_{p_*X}\hat\theta(p_*Y)=
\nabla^{\bar g}_X\bar\theta(Y)$.
Applying $\gamma^*$ to the both side of
\eqref{co},
we derive
\begin{equation*}
\begin{split}
&\gamma^*\bar g=\bar g=\mu\cdot p^*g.\\
&\gamma^*\mu\cdot\gamma^*p^*g=\gamma^*\mu\cdot p^*g.
\end{split}
\end{equation*}
Therefore 
$\gamma^*\mu=\mu$  which implies that
$\mu$ factors through a map $\hat \mu:M\ra \RR^+$
so that $p^*\hat g=p^*(\hat\mu\cdot g)$.
We have $\hat\mu\cdot g=\hat g$.
The conformal class of $g$ contains a l.c.K. metric $\hat g$
with parallel Lee form $\hat \theta$.
This finishes the proof of Theorem A. \hfill $\Box$

\vskip1cm

As to Corollary A$_1$ in the Introduction, we recall the following.
(Compare \cite{Va}, \cite[p.37]{DO}.)
Let $(M,g,J)$ be a compact, connected, non-K\"ahler, l.c.K. manifold
with parallel Lee form $\theta$.
Then the following results hold:
$g(\theta^\sharp,\theta^\sharp)=const$,
\begin{equation*}\begin{split}
&\Ll_{\theta^\sharp}J=\Ll_{J\theta^\sharp}J=0,\\
&\Ll_{\theta^\sharp}g=\Ll_{J\theta^\sharp}g=0.
\end{split}
\end{equation*}
Then $Z=\theta^\sharp-iJ{\theta^\sharp}$ is a holomorphic vector field
because $[\theta^\sharp,J\theta^\sharp]=0$ (cf. \cite{Ko1}).
By Definition \ref{parallellck},
$Z=\theta^\sharp-iJ{\theta^\sharp}$ is a holomorphic l.c.K. vector field.

\begin{pr}\label{real parallel flow}
The real vector fields $\theta^\sharp$ and $J{\theta^\sharp}$ satisfy
the following:
\begin{enumerate}
\item A flow generated by the Lee field $\theta^\sharp$
lifts to a one-parameter subgroup
of nontrivial homothetic holomorphic transformations w.r.t. $\Omega$.
\item A flow generated by the anti-Lee field $-J\theta^\sharp$
lifts to a one-parameter subgroup
consisting of holomorphic isometries w.r.t. $\Omega$.
\end{enumerate}
\end{pr}

\begin{proof}
Let $\{\hat\f_t\}_{t\in\RR}$ be the flow generated by $\theta^\sharp$ on $M$
and $\{\f_t\}_{t\in\RR}$ its lift to $\tilde M$. Denote by $\xi$
the vector field on $\tilde M$ induced by $\{\f_t\}$. Then,
$p_*\xi=\theta^\sharp$. Because $\theta$ is parallel, $\{\hat\f_t\}$ (resp.
$\{\f_t\}$) acts by holomorphic isometries w.r.t. $g$ (resp. $p^*g$).
In particular, $\{\f_t\}$ preserves $p^*\omega$. Then,
for $\Omega=e^{-\tau}p^*\omega$,
 we have $\displaystyle \f_t^*\Omega=e^{-(\f_t^*\tau-\tau)}\Omega$.
As $\rho:\{\f_t\}_{t\in\RR}\ra \RR^+$ is a homomorphism
and $\rho(\f_t)=e^{-(\f_t^*\tau-\tau)}$ is a constant for each $t\in\RR$
($\dim_\CC M\geq 2$), we can describe as
$-(\f_t^*\tau-\tau)=c\cdot t\ \ \mbox{for some constant}\ c.$
Recall that $h$ is the K\"ahler metric associated to $\Omega$. If $\{\f_t\}$ acts as holomorphic isometries w.r.t. $h$,
then the above equation implies that
$c=0$, $i.e.$ $\f_t^*\tau-\tau=0$ for every $t$, and so
$\mathcal{L}_\xi\tau=0$. On the other hand, as $d\tau=p^*\theta$, we have:
$$0=\mathcal{L}_\xi\tau=d\tau(\xi)=\theta(p_*\xi)=
\theta(\theta^\sharp)=const>0,$$
being a contradiction. Thus,
$\displaystyle \f_t^*\Omega=\rho(\f_t)\Omega =e^{c\cdot t}\Omega$
with $c\neq 0$.
Hence, $\{\f_t\}_{t\in \RR}$ is a group of nontrivial homothetic holomorphic
 transformations isomorphic to $\RR$.
On the other hand,
let $\{\hat\psi_t\}_{t\in\RR}$ (resp. $\{\psi_t\}_{t\in\RR}$)
 be the flow generated by -$J\theta^\sharp$ on $M$ (resp.
-$J\xi$ on $\tilde M$).
As $p_*(J\xi)=Jp_*\xi=J\theta^\sharp$,
$$
\mathcal{L}_{J\xi}\tau=d\tau(J\xi)=p^*\theta(J\xi)=\theta(J\theta^\sharp)=
g(J\theta^\sharp,\theta^\sharp)=0,
$$and hence $\psi_t^*\tau=\tau$ for every $t\in \RR$.
Using the fact that $\Ll_{J\theta^\sharp}g=0$,
$\Ll_{J\theta^\sharp}\omega=0$. This implies that
$\psi_t^*\Omega=\psi_t^*e^{-\tau}\psi_t^*p^*\omega=
e^{-\tau}p^*\hat\psi_t^*\omega=e^{-\tau}p^*\omega=\Omega.$
\end{proof}

Let $\RR\ra \tilde M\stackrel{\pi}{\lra}W$ be
the principal bundle where $\RR=\{\f_t\}_{t\in\RR}$ (cf. Lemma \ref{proper}).
Define the centralizer of $\RR$ in $\mathcal H(\tilde M,\Omega,J)$ to be:
\begin{de}\label{hol-omega}
$\mathcal C_{\mathcal H}(\RR)=\{f\in\mathcal H(\tilde M,\Omega,J)\ |\
f\circ \f_t=\f_t\circ f\ \ \mbox{for}\ \forall t\in \RR\}$.
\end{de}
As $\tilde{\mathcal C}$ centralizes the fundamental group $\pi_1(M)$,
 noting the remark below \eqref{assign},
\begin{equation}\label{pi-center}
\pi_1(M)\subset \mathcal C_{\mathcal H}(\RR).
\end{equation}

\begin{lm}\label{homomor}
There exists a homomorphism $\nu:\mathcal C_{\mathcal H}(\RR)
\ra{\rm PSH}(W,\eta,J)$ for which $\pi:\tilde M\ra W$ becomes
$\nu$-equivariant.
Moreover, there is a splitting homomorphism\\
$q: {\rm PSH}(W,\eta,J)\ra \mathcal C_{\mathcal H}(\RR)$.
\end{lm}

\begin{proof}
By definition, any element $f\in\mathcal C_{\mathcal H}(\RR)$
satisfies $f_*\xi=\xi$. As $f^*\Omega=\rho(f)\Omega$,
choosing $e^s=\rho(f)$, put $\gamma=\f_{-s}\circ f$. Then,
$\gamma^*\Omega=\Omega$. In particular, $\gamma$ leaves $W$ invariant.
Let $\gamma'$ be the restriction of $\gamma$ to $W$
$(i.e.,\ i\circ \gamma'=\gamma)$. Using \eqref{et} and $\gamma_*\xi=\xi$,
we have that
${\gamma'}^*\eta=\gamma^*\Ll_\xi\Omega=\Ll_\xi\Omega=\eta.$
Hence $\gamma'\in {\rm PSH}(W,\eta,J)$.
If we define $\nu(f)=\gamma'$, then it is easy to see that
$\nu$ is a well defined homomorphism.
Let $x=\f_tw$ be a point in $\tilde M$. As $\pi(x)=w$,
$\pi(fx)=\pi(\f_s\gamma(\f_tw))=
\pi(\f_s\f_t i\gamma'w)=\pi(i\gamma'w)=\gamma'w=\nu(f)\pi(x)$, so $\pi$
is $\nu$-equivariant.

For $\gamma\in{\rm PSH}(W,\eta,J)$, we define a diffeomorphism 
$\tilde \gamma:\tilde M\ra \tilde M$ to be
\begin{equation}\label{split-action}
\tilde \gamma(x)=\tilde \gamma(\f_tw)=\f_t \gamma w.
\end{equation}By definition, $\pi\circ \tilde \gamma=\gamma\circ \pi$ and
the $t$-coordinate satisfies that $\tilde \gamma^*t=t$.
Using \eqref{eta} and $\gamma^*\eta=\eta$, it follows that
$\tilde \gamma^*\Omega=d(e^{\gamma^*t}\pi^*\gamma^*\eta)=
d(e^t\pi^*\eta)=\Omega.$
To see that $\tilde \gamma:\tilde M\ra \tilde M$ is holomorphic,
notice that $\tilde \gamma_*\xi=\xi$. As
$\tilde \gamma(\psi_\theta x)= \tilde \gamma(\psi_\theta\f_tw)
=\tilde \gamma(\f_ti{\psi'}_\theta w)=\f_t i\gamma{\psi'}_\theta w$,
and ${\gamma}_*A=A$,
\begin{equation}\label{diff1}
\begin{split}
&\tilde \gamma_*(-J\xi_x)=\tilde \gamma_*(\frac{d\psi_\theta}{d\theta}(x)|_{\theta=0})=
(\frac{d\f_t i\gamma({\psi'}_\theta w)}{d\theta}|_{\theta=0})\\
&={\f_t}_*i_*\gamma_*(\frac{d{\psi'}_\theta}{d\theta}(w)|_{\theta=0})
={\f_t}_*i_*\gamma_*A_w={\f_t}_*i_*A_{\gamma w}={\f_t}_*(-J\xi_{\gamma w})=
-J\xi_{\tilde \gamma x}.
\end{split}
\end{equation}
Hence, $\tilde \gamma$ preserves $\{\xi,J\xi\}^{\perp}$.
Since the complex structure $J:\mathop{\Null}\, \eta\ra {\Null}\, \eta$
is defined by the
commutative diagram \eqref{J-com},
$J\gamma_*(\pi_*X)=\gamma_*J(\pi_*X)$ for $X\in\{\xi,J\xi\}^{\perp}$
 by definition.
Then $\pi_*\tilde \gamma_*J(X)=
J\gamma_*\pi_*(X)=J\pi_*\tilde \gamma_*(X)=\pi_*J\tilde \gamma_*(X)$.
As a consequence, $\tilde \gamma_*\circ J=J\circ \tilde \gamma_*$
 on $\tilde M$.
Hence, $\tilde \gamma\in\mathcal C_{\mathcal H}(\RR)$.
It is easy to check that $q(\gamma)=\tilde \gamma$ is a
homomorphism of 
${\rm PSH}(W,\eta,J)$ into $\mathcal C_{\mathcal H}(\RR)$
such that $\nu\circ q={\rm id}$.

\end{proof}

\begin{re}\label{Seifertaction}\ \ From this lemma,
there is an isomorphism $\mathcal C_{\mathcal H}(\RR)
\approx\RR\times {\rm PSH}(W,\eta,J)$ where
each element of $\mathcal C_{\mathcal H}(\RR)$
is described as $\f_s\cdot q(\al)$ for $s\in \RR,\ \al\in{\rm PSH}(W,\eta,J)$.
It acts on $\tilde M$ as
\[
\f_s\cdot q(\al)(\f_t\cdot w)=
\f_{s+t}\cdot \al w,
\]for which there is an equivariant principal bundle:
\[
\RR\ra (\mathcal C_{\mathcal H}(\RR),\tilde M)\stackrel
{(\nu,\pi)}\lra  ({\rm PSH}(W,\eta,J),W).
\]
\end{re}

\subsection{Central group extension}\label{Central groupex}
Consider the exact sequence:
\begin{equation}\label{exact}
1\ra \RR\ra \mathcal C_{\mathcal H}(\RR)\stackrel{\nu}{\lra}
{\rm PSH}(W,\eta,J)\ra 1.
\end{equation}
Suppose that $\RR\cap \pi_1(M)$ is nontrivial. Then it is
an infinite cyclic subgroup $\ZZ$ such that
the quotient group $\RR/\ZZ$ is a circle $S^1$.
Put $Q=\nu(\pi_1(M))\subset {\rm PSH}(W,\eta,J)$.
We have a central group extension:
\begin{equation}\label{central}
1\ra\ZZ\ra \pi_1(M)\stackrel{\nu}{\lra}Q\ra 1.
\end{equation}The above principal bundle restricts to
the following one:
\begin{equation}\label{group extension}
(\ZZ,\RR)\ra (\pi_1(M), \tilde M)\stackrel{(\nu,\pi)}{\lra}
(Q,W).
\end{equation}
As both $\RR$ and $\pi_1(M)$ act properly on $\tilde M$,
$Q$ acts also properly discontinuously
(but not necessarily freely) on $W$ such that
the quotient Hausdorff space $W/Q$ is compact.
Since $\rho(\ZZ)\subset \rho(\RR)=\RR^+$ from $\S$ \ref{sub-ps},
$\rho(\ZZ)$ is an infinite cyclic subgroup of $\RR^+$.
We need the following lemma. (Compare \cite{kamicomp}, \cite{CR}.)

\begin{lm}\label{splits}
Let $1\ra\ZZ\ra \pi_1(M)\stackrel{\nu}{\lra}Q\ra 1$
be the central extension as in \eqref{central}.
Then, $\pi_1(M)$ has a splitting subgroup $\pi'$ of finite index:
$1\ra\ZZ\ra \pi'\stackrel{\nu}{\lra}Q'\ra 1$
In particular, there exists a subgroup $H'$ of $\pi'$
which maps isomorphically onto a subgroup $Q'$ of
finite index in $Q$.
\end{lm}

\begin{proof} 
Consider the homomorphism $\rho'=\rho|_{\pi_1(M)}: \pi_1(M)\lra \RR^+$
from \eqref{assign}. Then, $\rho'(\pi_1(M))$ is a free abelian group of
rank $k\geq 1$. If we note that $\rho'(\ZZ)$ is an infinite cyclic subgroup
of $\rho'(\pi_1(M))$, then we can choose a subgroup $G$ of finite index
in $\rho'(\pi_1(M))$ such that $\rho'(\ZZ)$
is a direct summand in $G$;\ $G=\rho'(\ZZ)\times \ZZ^{k-1}$. 
Put $\pi'={\rho'}^{-1}(G)$ and $H'={\rho'}^{-1}(\ZZ^{k-1})$.
Then, $\pi'$ has finite index in $\pi_1(M)$. Obviously
$\nu$ maps $H'$ isomorphically onto $\nu(H')=Q'$
which is of finite index in $Q$.

\end{proof}

\begin{pr}\label{smoothcr}
The subgroup $Q'$ acts freely on $W$ so that
the orbit space $W/Q'$ is a closed strictly pseudoconvex 
pseudo-Hermitian manifold induced from the\\
pseudo-Hermitian structure $(\eta,J)$ on $W$.
\end{pr} 

\begin{proof}
Let $f={\nu'}^{-1}: Q'\ra H'$ be the inverse isomorphism.
For each $\al'\in Q'$ there exists a unique element $\la(\al')\in \RR$
such that $f(\al')=\f_{\la(\al')}\cdot q(\al')$.
As we know that $Q$ acts properly discontinuously
on $W$ from the remark below \eqref{group extension},
the stabilizer at each point is finite.
Suppose that $\al'w=w$ for some point $w\in W$.
As $\al'\in Q_w$, ${\al'}^\ell=1$ for some $\ell$.
Since $\f_t$ is a central element and $q$ is a homomorphism,
$1=f({\al'}^\ell)=\f_{\ell\la(\al')}\cdot q({\al'}^{\ell})
=\f_{\ell\la(\al')}$. Thus, $\la(\al')=0$, $i.e.,$ $f(\al')=q(\al')$.
By definition of the action $(\pi',\tilde M)$,
$f(\al')(\f_tw)=q(\al')(\f_tw)=\f_t \al' w=
\f_t w$. As $\pi'$ acts freely on $\tilde M$,
$f(\al')=1$ and so $\al'=1$.
If we note that $Q'\subset{\rm PSH}(W,\eta,J)$,
then $(\eta,J)$ induces a pseudo-Hermitian structure $(\hat\eta,J)$ on $W/Q'$.
Here we use the same notation $J$ to the complex structure on 
$\mathop{\Null}\ \hat \eta$.

\end{proof}
 
\section{Examples of l.c.K. manifolds with parallel Lee form}
In this section we present an explicit construction for the Hopf manifolds.\\
Let  $S^{2n-1}=\{(z_1,\ldots ,z_n)\in \CC^n\ |\  \abs{z_1}^2+\cdots+\abs{z_n}^2=1\}$ be the 
sphere endowed with its  standard contact structure
\begin{equation}\label{stand}
\eta_0=\sum_{j=1}^n(x_jdy_j-y_jdx_j), \; \text{where}\; z_j=x_j+\sqrt{-1}\ y_j.
\end{equation} 
Let $J_0$ be the restriction of the standard complex structure of 
$\CC^n$ to $\CC^n-\{0\}$. It is known that
 the group of pseudo-Hermitian transformations, 
$\mathrm{PSH}(S^{2n-1},\eta_0,J_0)$
is isomorphic with $\mathrm{U}(n)$ (see \cite{W1}, for example). 
We define a $1$-parameter subgroup\\
 $\{\psi_t\}_{t\in \RR}\subset \mathrm{PSH}(S^{2n-1},\eta_0,
J_0)$ by the formula:
$$\psi_t(z_1,\ldots, z_n)=(e^{\mathrm{i}ta_1}z_1,\ldots,  e^{\mathrm{i}ta_n}z_n),$$
where $\mathrm{i}=\sqrt{-1}$ and $a_1,\ldots,a_n\in\RR$. The vector field induced by this 
action is 
$$A=\sum_{j=1}^na_j(x_j\frac{d}{dy_j}-y_j\frac{d}{dx_j})$$
and satisfies $\eta_0(A)=a_1\abs{z_1}^2+\cdots+a_n\abs{z_n}^2.$

Now we require that $\eta_0(A)>0$ everywhere on $S^{2n-1}$.
 Then the numbers $a_k$ must satisfy (up to rearrangement):
\begin{equation}\label{ups}
0<a_1\leq\cdots\leq a_n.
\end{equation}
Define a new contact form $\eta_A$ on the sphere by 
$$\eta_A=\frac{1}{{\sum}_{j=1}^n a_j\abs{z_j}^2}\cdot\eta_0.$$
The contact distributions of $\eta_0$ and $\eta_A$ coincide, but the characteristic field of 
$\eta_A$ is $A$: $\eta_A(A)=1$, $\iota_Ad\eta_A=0$. As $A$ generates the flow $\{\psi_t\}_{t\in\RR}\subset\mathrm{PSH}(S^{2n-1},\eta_0,J_0)$,
 note that
$\psi_{t*}\circ J_0=J_0\circ \psi_{t*}$ on $\Null{\eta_A}$. 
Define a $2$-form on the product $\RR\times S^{2n-1}$ by:
\[
\Omega_A=2d(e^t{\rm pr}^*\eta_A),\ \ (t\in \RR).
\]Here ${\rm pr}:\RR\times S^{2n-1}\ra S^{2n-1}$ is the projection.
If $\RR=\{\f_s\}_{s\in\RR}$ acts on $\RR\times S^{2n-1}$ by left
translations: $\f_s(t,z)=(s+t,w)$,
 then the group $\RR\times \mathrm{PSH}(S^{2n-1},\eta_A,J_0)$ acts by 
homothetic transformations w.r.t. $\Omega_A$:
\begin{equation}\label{homothetic}
(\f_s\times \al)^*\Omega_A=e^s\cdot \Omega_A, \quad (\al\in 
\mathrm{PSH}(S^{2n-1},\eta_A,J_0)).
\end{equation} 
In general,  $\mathrm{PSH}(S^{2n-1},\eta_A,J_0)$
is the centralizer of $\{\psi_t\}_{t\in\RR}$ in $\mathrm{U}(n)$.
In view of the formula of $\psi_t$, $\mathrm{PSH}(S^{2n-1},\eta_A,J_0)$
contains the maximal torus of $\mathrm{U}(n)$ at least.
\begin{equation}\label{maximaltorus}
T^n\subset \mathrm{PSH}(S^{2n-1},\eta_A,J_0).
\end{equation}
(For example, if all $a_j$ are distinct, 
$\mathrm{PSH}(S^{2n-1},\eta_0,J_0)=T^n$).

Let $\displaystyle N=\frac{d}{dt}$ be the vector field induced on 
$\RR\times S^{2n-1}$  by the $\RR$-action.
Taking into account that
$T(\RR\times S^{2n-1})=N\oplus A\oplus \Null{\eta_A}$, we define
an almost complex structure $J_A$ on $\RR\times S^{2n-1}$ by:
\begin{equation*}
\begin{split}
&J_AN=-A, \quad J_AA = N,\\
&\displaystyle J_A|{\Null{\eta_A}}=J_0
\end{split}
\end{equation*} 
and show its integrability. Indeed, let
\[
T(\RR\times S^{2n-1})\otimes\CC=
\{T^{1,0}+(A-\mathrm{i}N)\}\oplus \{T^{0,1}+(A+\mathrm{i}N)\}
\]be the splitting corresponding 
to $J_A$ (here $T^{1,0}+T^{0,1}=\Null{\eta_A}\otimes\CC$).
As $\displaystyle J_A|{\Null{\eta_A}}=J_0$, $[T^{1,0},T^{0,1}]\subset T^{1,0}$. Recalling that $A$ is the characteristic field of $\eta_A$,
we see that\\
$[X,A]\in \Null{\eta_A}$ for any $X\in\Null{\eta_A}$.
If $X\in T^{1,0}$, 
then $[X,A-\mathrm{i}N]=[X,A]=\displaystyle 
\lim_{t\rightarrow 0}\frac{X-\psi_{-t*}X}{t}$.
Noting that $\psi_t\in\mathrm{PSH}(S^{2n-1},\eta_A,J_0)$ $(i.e.,\
\psi_{t*}J_0=J_0\psi_{t*})$, 
\begin{equation*}
\begin{split}
J_A[X,A-\mathrm{i}N]&=J_0[X,A]=\lim_{t\rightarrow 0}\frac{J_0X-\psi_{-t*}J_0X}{t}=[J_0X,A]\\
&=[\mathrm{i}X,A]=\mathrm{i}[X,A]=\mathrm{i}[X,A-\mathrm{i}N].
\end{split}
\end{equation*}
Thus $[X,A-\mathrm{i}N]\in \{T^{1,0}+(A-\mathrm{i}N)\}$.
Hence $J_A$ is integrable.
By the definition of $J_A$,
it is easy to check that
the elements of $\RR\times\mathrm{PSH}(S^{2n-1},\eta_A,J_0)$ 
are holomorphic w.r.t. $J_A$.
Moreover, $\Omega_A$ is $J_A$-invariant. Hence,
$\Omega_A$ is a K\"ahler form on the complex manifold
$(\RR\times S^{2n-1},J_A)$ on which
$\RR\times\mathrm{PSH}(S^{2n-1},\eta_A,J_0)$ acts as the group of 
holomorphic homothetic transformations. Define a Hermitian metric
$\tilde g_A$ and its
fundamental $2$-form $\tilde \omega_A$ by setting
\begin{equation}\label{lckpara}
\begin{split}
&\tilde\omega_A=2e^{-t}\cdot\Omega_A.\\
&\tilde g_A(X,Y)=\tilde \omega_A(J_AX,Y), \quad 
\forall\ X,Y\in T(\RR\times S^{2n-1}).
\end{split}
\end{equation}(Compare \eqref{bar}.)
By \eqref{homothetic},
$\RR\times\mathrm{PSH}(S^{2n-1},\eta_A,J_0)$ acts as holomorphic 
isometries of $(\tilde g_A,J_A)$.
When we choose a properly
discontinuous group $\Gamma\subset \RR\times\mathrm{PSH}(S^{2n-1},\eta_A,J_0)$ acting freely on $\RR\times S^{2n-1}$, $\tilde g_A$ (resp. $\tilde \omega_A$)
 induces
a Hermitian metric $g_A$ (resp. the fundamental $2$-form $\omega_A$)
on the quotient complex manifold
$(\RR\times S^{2n-1}/\Gamma,\hat J_A)$, 
where the complex structure $\hat J_A$ is induced from $J_A$. 
We have to check that $g_A$ is a l.c.K. metric with parallel Lee form.
Let $p:\RR\times S^{2n-1}\ra\RR\times S^{2n-1}/\Gamma$ be the
projection so that $p^*\omega_A=\tilde \omega_A$.
Since $\tilde\omega_A=e^{-t}\cdot\Omega_A$,
we have $d\tilde\omega_A=-dt\wedge\tilde \omega_A$. Thus  
$\tilde g_A$ is a l.c.K. metric 
with Lee form $d(-t)$ on $\RR\times S^{2n-1}$.
If we note that the group 
$\RR\times\mathrm{PSH}(S^{2n-1},\eta_A,J_0)$
leaves $d(-t)$ invariant, $i.e.\ (\f_s\times \al)^*d(-t)=d(-(s+t))=d(-t)$,
then $d(-t)$ induces a $1$-form $\theta$ on 
$\RR\times S^{2n-1}/\Gamma$ such that $p^*\theta=d(-t)$. 
The equation $d\tilde\omega_A=-dt\wedge\tilde \omega_A$ implies that
$d\omega_A=\theta\wedge\omega_A$ on
$\RR\times S^{2n-1}/\Gamma$. As $d\theta=0$,
$g_A$ is a l.c.K. metric with Lee form $\theta$.
For the rest, the same argument as in the proof of
Lemma \ref{parallelform} can be applied to show that
$\theta$ is the parallel Lee form of $g_A$.
Finally, we examine the complex structure
$\hat J_A$ on $\RR\times S^{2n-1}/\Gamma$.
Let  $H:\RR\times S^{2n-1}\rightarrow \CC^n-\{0\}$ be the 
diffeomorphism defined by:
$$H(t, (z_1,\ldots,z_n))=(e^{-a_1t}z_1,\ldots,  e^{-a_nt}z_n),$$
where $\{a_1,\ldots,a_n\}$ satisfies the condition \eqref{ups}. We shall show that $H$ is a 
$(J_A,J_0)$-biholomorphism. We have:
\begin{equation*}
\begin{split}
H_*(N_{(s,z)})&=\frac{dH({t+s},z)}{dt}{|_{t=0}}
=(-a_1\cdot e^{-a_1s}\cdot z_1,\ldots,-a_n\cdot e^{-a_ns}\cdot z_n);\\
H_*(J_AN_{(s,z)})&
=H_*(-A_{(s,z)})=-H_*((s,\frac{d}{dt}(e^{\mathrm{i}ta_1}z_1,\ldots,
e^{\mathrm{i}ta_n}z_n)|_{t=0})\\
&=-(\mathrm{i}a_1e^{-a_1s}z_1,\ldots,\mathrm{i}a_ne^{-a_ns}z_n)=
J_0H_*(N_{(s,z)}).
\end{split}
\end{equation*}
From $H_*(A_{(s,z)})=-J_0H_*(N_{(s,z)})$, we derive 
$J_0H_*(A_{(s,z)})=H_*(N_{(s,z)})=H_*(J_AA).$
Now let $X\in \Null{\eta_A}\subset TS^{2n-1}$ and let $\sigma(t)$ be an 
integral curve of $X$ on $S^{2n-1}$: $\dot\sigma(t)=X$, $\dot \sigma(0)=X_z$. 
We can view $X$ as a pair: $X_{(s,z)}=(s,\dot\sigma(0))$. Then:
$$H_*(X_{(s,z)})=\frac{d}{dt}H(s, \sigma(t))|_{t=0}=
(e^{-a_1s}\dot\sigma_1(0),\ldots,e^{-a_ns}\dot\sigma_n(0)).$$
From this we obtain:
\begin{equation*}
\begin{split}
H_*(J_AX_{(s,z)})&=H_*((s,J_0\dot\sigma(0)))=H_*((s,
(\mathrm{i}\dot\sigma_1(0),\ldots,\mathrm{i}\dot\sigma_n(0))))\\
&=(\mathrm{i}e^{-a_1s}\dot\sigma_1(0),\ldots,\mathrm{i}e^{-a_ns}\dot\sigma_n(0))\\&=J_0(e^{-a_1s}
\dot\sigma_1(0),\ldots,e^{-a_ns}\dot\sigma_n(0))
=J_0H_*(X_{(s,z)}).
\end{split}
\end{equation*}
Therefore $H:(\RR\times S^{2n-1},J_A)\rightarrow (\CC^n-\{0\},J_0)$ is a
 biholomorphism.

Let $\mathrm{Hol}(\CC^n-\{0\},J_0)$ be the group of 
all biholomorphic transformations. If we associate to each
 $\gamma\in \RR\times\mathrm{PSH}(S^{2n-1},\eta_A,J_0)$ the biholomorphic map
 $H\circ\gamma\circ H^{-1}$, we obtain a faithful
homomorphism
$\RR\times\mathrm{PSH}(S^{2n-1},\eta_A,J_0)
\lra \mathrm{Hol}(\CC^n-\{0\},J_0)$. Let $\Gamma_H$ be
the image of $\Gamma$ in $\mathrm{Hol}(\CC^n-\{0\},J_0)$.
\begin{de}
The quotient complex manifold $\CC^n-\{0\}/\Gamma_H$ is called 
a Hopf manifold.
\end{de}
\noindent We have shown:
\begin{te}\label{secondary}
The Hopf manifold
$\CC^n-\{0\}/\Gamma_H$ admits a l.c.K. metric $g$
with parallel Lee form $\theta$.
\end{te}

By \eqref{maximaltorus},
$T^n\subset \mathrm{PSH}(S^{2n-1},\eta_A,J_0)$.
Choose $s\in \RR-\{0\}$ and $n$-complex numbers $c_1,\ldots,c_n\in S^1$.
Consider an infinite cyclic subgroup $\ZZ$
generated by the element\\
$(s,(c_1,\ldots,c_n))$ from 
$\RR\times \mathrm{PSH}(S^{2n-1},\eta_0,J_0)$.
Then the corresponding group 
$\ZZ_H$ is generated by the element
$(e^{-a_1s}\cdot c_1,\ldots,e^{-a_ns}\cdot c_n)$
acting on $\CC^n-\{0\}$.
Let $\Lambda=(\la_1,\ldots,\la_n)$, with $\la_j=e^{-a_js}\cdot c_j$
and so $\ZZ_H=<(\la_1,\ldots,\la_n)>$.
The condition \eqref{ups} ensures that the complex numbers 
$\la_j$ satisfy 
$$0<\abs{\la_n}\leq\cdots\leq\abs{\la_1}<1.$$
Put $M_\Lambda=\CC^n-\{0\}/\Gamma_H$.
We call $M_\Lambda$ a \emph{primary Hopf manifold of type $\Lambda$}. Indeed, for $n=2$, one recovers the 
primary Hopf surfaces of K\"ahler rank $1$. 
In particular, we derive Theorem B in the Introduction.

\begin{re}
Note that the manifolds  $M_\Lambda$ are all diffeomorphic with $S^1\times S^{2n-1}$ and that for $c_1=\cdots=c_n=1$ and $a_1=\cdots =a_n$, 
we obtain the standard Hopf manifold, the first known example of a l.c.K. manifold with parallel Lee form, cf. \cite{Va}.\\
In \cite{GO} a l.c.K. metric with parallel Lee form is constructed on the primary Hopf surface $M_{\la_1,\la_2}=\CC^2-\{0\}/\Gamma$, $\Gamma\cong \ZZ$ generated by $(z_1,z_2)\mapsto (\la_1z_1,\la_2z_2)$, $\abs{\la_1}\geq\abs{\la_2}>1$. There the diffeomorphism between $M_{\la_1,\la_2}$ and $S^1\times S^3$ is used to construct a potential for the K\"ahler metric $h$ (in the present paper notations) on the universal cover. The same  diffeomorphism is then used to transport  the l.c.K. structure on $S^1\times S^3$ and to show that the  induced Sasakian structure on $S^3$ is  a deformation of the standard Sasakian structure of the $3$-sphere. See also \cite{Be} where a complete list of compact, complex surfaces admitting l.c.K. metrics with parallel Lee form is provided.
\end{re} 

\section{Lee-Cauchy-Riemann transformations}
In this section, we  consider the group $\mathop{\Aut}_{LCR}(M)$ described in 
the Introduction.\\
Let  $\{\theta, \theta\circ J, \theta^\al,\bar\theta^\al\}_{\al=1,
\cdots, n-1}$ be a unitary,
local coframe field 
adapted to a l.c.K. manifold $(M,g,J)$.
Consider the subgroup $G$ of $\mathop{\GL}(2n,\RR)$
consisting of the following elements:

$$
\left\{\left(\begin{array}{cccccc}
1&0&\ & 0        &\  &0\\
0&u&\ &v^{\alpha}&\ &\bar v^{\alpha}\\
0&0&\ &\sqrt u \, U_{\beta}^{\alpha}&\ &0\\
0&0&\ &0&\ &\sqrt u \, \bar U_{\beta}^{\alpha}
\end{array}\right)
\ \mid \ u\in\RR^+,v^{\alpha}\in\CC, U_{\beta}^{\alpha}\in {\rm 
U}(n-1)\right\}. 
$$
Let $G\ra P\ra M$  be the principal bundle of
the $G$-structure consisting of the above coframes
$\{\theta, \theta\circ J, \theta^{\alpha},\bar\theta^{\alpha}\}$.
If we note that $G$ is isomorphic to the semidirect
product\\
$\CC^{n-1}\rtimes ({\rm U}(n-1)\times \RR^+)$,
then the Lie algebra ${\mathfrak g}$
is isomorphic to $\CC^{n-1}+{\mathfrak u}(n-1)+\RR$.
In particular, the matrix group
${\mathfrak g}\subset {\mathfrak gl}(2n,\RR)$
has no element of ${\rm rank}$
 1, \emph{i.e.} it is  \emph{elliptic} (cf. \cite{Ko}).
Note that
$\CC^{n-1}$ is of infinite type, while ${\mathfrak u}(n-1)+\RR$
is of order $2$. As $M$ is assumed to be
compact, the group of automorphisms $\mathcal U$ of
$P$ is a (finite dimensional) Lie group.
\begin{de}
The group of all diffeomorphisms of $M$ onto itself
which preserve the above $G$-structure is denoted by
$\mathop{\Aut}_{LCR}(M,g,J,\theta)$ 
(or simply by $\mathop{\Aut}_{LCR}(M)$). 
We call $\mathop{\Aut}_{LCR}(M)$
the group of Lee-Cauchy-Riemann transformations
on a l.c.K. manifold $(M,g,J)$ adapted to the Lee form $\theta$.
\end{de}
By  definition, if $f\in\mathop{\Aut}_{LCR}(M)$, then $f^* :P\ra P$
is a bundle automorphism satisfying 
\begin{equation} \label{rel}
\begin{split}
&f^*\theta=\theta,\\
&f^*(\theta\circ J)=\lambda\cdot(\theta\circ J),\ \mbox{for some positive,
smooth function}\ \lambda,\\
&f^*\theta^{\al}=\sqrt \lambda\cdot\theta^\be V^{\al}_{\be}+
(\theta\circ J)\cdot w^{\al},\\
&f^*{\bar \theta}^{\al}=\sqrt \lambda\cdot
{\bar\theta}^\be{\bar V}^\al_\be+(\theta\circ J)\cdot{\bar w}^{\al},
\end{split}
\end{equation}
for functions $V^\al_\be\in {\rm U}(n-1)$ and 
$w^\al\in \CC$.
Note that the group of holomorphic isometries
${\rm I}(M,g,J)$ is contained in $\mathop{\Aut}_{LCR}(M)$.
In fact, an element $f\in {\rm I}(M,g,J)$
satisfies  $f^*\theta=\theta$, $f^*(\theta\circ J)=(\theta\circ J)$
and $f^*\omega=\omega$. Let  $\{\theta^{\sharp}, J\theta^{\sharp}\}^{\perp}$
be the orthogonal complement of
the complex plane field $\{\theta^{\sharp}, J\theta^{\sharp}\}$ w.r.t. $g$.
It is obviously $J$-invariant. If we note that $\omega|\{\theta^{\sharp}, 
J\theta^{\sharp}\}^{\perp}=-\mathrm{i}\sum_{\al,\be}\delta_{\al\be}\theta^\al\wedge\bar\theta^\be$, then 
$f^*\theta^{\al}=\theta^\be U^{\al}_{\be}$, $f^*{\bar \theta}^{\al}={\bar\theta}^\be{\bar U}^\al_\be$
for some  matrix function $U^\al_\be\in {\rm U}(n-1)$. 

\begin{lm}\label{perp}
Any element $f\in\mathop{\Aut}_{LCR}(M)$  
preserves $\{\theta^{\sharp}, J\theta^{\sharp}\}^{\perp}$ and
is holomorphic on it.
\end{lm}
\begin{proof}
Let $X\in\{\theta^{\sharp}, J\theta^{\sharp}\}^{\perp}$.
The equations $f^*\theta=\theta$, $f^*(\theta\circ J)=\lambda\cdot
(\theta\circ J)$ show that
\begin{equation}\label{orthonomal}
\begin{split}
 g(f_*X,\theta^{\sharp})&=\theta(f_*X)=\theta(X)=g(X,\theta^{\sharp})=0,\\
g(f_*X,J\theta^{\sharp})&=-g(Jf_*X,\theta^{\sharp})=-\theta(Jf_*X)
=-\theta\circ J(f_*X)\\
&=-\lambda\cdot \theta\circ J(X)=-g(X,(\theta\circ J)^{\sharp})
=g(X,J\theta^{\sharp})=0.
\end{split}
\end{equation}
Thus $f_*$ applies $\{\theta^{\sharp}, J\theta^{\sharp}\}^{\perp}$ onto itself. Moreover, if $\displaystyle \theta_\al^{\sharp}$ is a dual frame field to $\theta^\al$
(similarly for $\bar\theta^{\al}$), then
the frame $\displaystyle \{\theta_\al^{\sharp},\bar\theta_{\al}^{\sharp}\}_{\al=1,\cdots,n-1}$
spans $\{\theta^{\sharp}, J\theta^{\sharp}\}^{\perp}\otimes\CC$.\\
The equation $f^*\theta^{\al}=\sqrt \lambda\cdot\theta^\be V^{\al}_{\be}+
(\theta\circ J)\cdot w^{\al}$ implies that
$f_*\theta_{\al}^{\sharp}=\sqrt \lambda\cdot\theta_\be^{\sharp} V_{\al}^{\be}$
\ (similary for $f_*\bar\theta_{\al}^\sharp$). Therefore $f_*\circ J=J\circ f_*$on $\{\theta^{\sharp}, J\theta^{\sharp}\}^{\perp}$.

\end{proof}

When  a noncompact $LCR$ flow exists on a compact l.c.K. manifold $M$
with parallel Lee form,
we shall prove a rigidity similar 
to the  one implied by a noncompact $CR$-flow
on a compact $CR$-manifold (cf. \cite{Ob}, \cite{Kam}). 

\subsection*{Proof of Theorem C}\hfill 
\subsection{Existence of spherical $CR$-structure
on $W/Q'$} 
Let $1\ra\ZZ\ra \pi'\stackrel{\nu}{\lra}Q'\ra 1$
be the split central group extension from Lemma \ref{splits}.
Put $M'=\tilde M/\pi'$. Then it is easy to see that the Lee form $\theta$,
the LCR-action $\CC^*$ lift to those of $M'$,
 so we retain the same notations for $M'$.
We put $\CC^*=S^1\times \RR^+$ where
$\RR^+=\{\hat\h_t\}_{t\in\RR}$ is a $LCR$ flow on $M'$.  By hypothesis,
$S^1=\{\hat\varphi_t\}_{t\in\RR}$ induces the Lee field $\theta^{\sharp}$.
From 1 of Proposition \ref{real parallel flow},
 $S^1$ lifts to a nontrivial holomorphic homothetic flow
$\RR=\{\varphi_t\}_{t\in\RR}$ on $\tilde M$ w.r.t. $\Omega$.
We obtain a LCR-action of $\RR\times \RR^+$ on $\tilde M$
for which $\RR$ acts properly as before.
Consider the commutative diagram of principal bundles:
\begin{equation}\label{comm}
\begin{CD}
\mathbb{Z}@>>>\pi'@>\nu>>Q'\\
@VVV          @VVV        @VVV\\
\RR@>>>(\RR\times \RR^+, \tilde M)@>(\tilde \nu,\pi) >>(\RR^+, W)\\
@VVV @VVpV @VVpV\\
S^1 @>>> (S^1\times \RR^+, M') @>(\hat\nu,\hat\pi) >> (\RR^+, W/Q')
\end{CD}
\end{equation}From the bottom line, the projection $\hat\nu$
maps the group $\RR^+=\{\hat \h_{t}\}_{t\in\RR}$ onto 
a group $\RR^+=\{\bar\h_t\}_{t\in\RR}$ acting on $W/Q'$.

\begin{lm}\label{cr-trans}
The group $\RR^+=\{\bar\h_t\}_{t\in\RR}$ acts by $CR$-transformations
on $W/Q'$ w.r.t. the $CR$-structure induced from the
strictly pseudoconvex, pseudo-Hermitian structure
$(\hat\eta,J)$.
\end{lm}

\begin{proof}
As $\xi$ generates the flow
$\RR=\{\varphi_t\}_{t\in\RR}$,
$p_*\xi=\theta^{\sharp}$ on $M'$ by hypothesis and so
$p:\tilde M\ra M'$
maps the complex plane field $\{\xi,J\xi\}$ onto
$\{\theta^{\sharp},J\theta^{\sharp}\}$.
By Lemma \ref{perp},
each $\hat\h_t\in\mathop{\Aut}_{LCR}(M')$
preserves $\{\theta^\sharp, (\theta\circ J)^\sharp\}^\perp$.
So its lift $\h_t$ preserves
the $J$-invariant distribution
$\{\xi, J\xi\}^\perp$.
Since $\pi_*:(\{\xi, J\xi\}^\perp,J)\ra (\mathop{\Null}\, \eta,J)$
is $J$-isomorphic and
 each $\h_t$ is holomorphic on $\{\xi, J\xi\}^\perp$,
$\hat\pi_*:(\{\theta^\sharp, (\theta\circ J)^\sharp\}^\perp,J)
\lra (\mathop{\Null}\, \hat \eta,J)$
is also $J$-isomorphic through the commutative diagram
and thus each $\bar\h_t$ is holomorphic on
$\mathop{\Null\ \hat\eta}$; $(\bar\h_{t*}\circ J= J\circ \bar\h_{t*})$.
Therefore, $\RR^+=\{\bar\h_t\}_{t\in\RR}$ is a closed, noncompact
subgroup of $CR$-transformations of $W/Q'$ w.r.t.
$(\mathop{\Null}\, \hat\eta,J)$.

\end{proof}

By this lemma, we obtain a
compact strictly pseudoconvex
$CR$-manifold $W/Q'$ admitting 
a closed, noncompact $CR$-transformations 
$\RR^+$. Then we apply the result of \cite{Kam}
to show that $W/Q'$ is $CR$-equivalent to the sphere
$S^{2n-1}$ with the standard $CR$-structure.
In particular $Q'=\{1\}$ and thus $Q$ is a finite subgroup of 
${\rm PSH}(W,\eta,J)$ from Lemma \ref{splits}.
By definition of spherical $CR$-structure (cf. \cite{Ku}, \cite{Gold}), 
there exists a developing pair:
\[
(\mu,\, \mathop{\dev}):
({\Aut}_{CR}(W),W)\ra ({\rm PU}(n,1),S^{2n-1})
\]for which $\mathop{\dev}$ is a $CR$-diffeomorphism and
$\mu:\mathop{\Aut}_{CR}(W)\ra {\rm PU}(n,1)$
is the holonomy isomorphism.
Here $\mbox{\rm PU}(n,1)=\mathop{\Aut}_{CR}(S^{2n-1})$
and $\mathop{\Aut}_{CR}(W)$  
is the group of all $CR$-automorphisms of $W$
containing the groups $\RR^+$ and ${\rm PSH}(W,\eta,J)\supset Q$.\\
As $S^1\ (\subset \CC^*)$ acts on $M$ without fixed points
(but not necessarily freely),
the quotient space $M/S^1=W/Q(\approx S^{2n-1}/\mu(Q))$
is an orbifold, so such a finite subgroup $Q$ may exist.

On the other hand, we recall some facts from the theory of hyperbolic groups
(cf. \cite{CG}).
The noncompact closed $\mu(\RR^+)$-action on $S^{2n-1}$ is characterized as
whether it is either loxodromic $(=\RR^+)$
or parabolic $(=\mathcal R)$ for which
$\RR^+$ has exactly two fixed points $\{0,\infty\}$ or
$\mathcal R$ has the unique fixed point $\{\infty\}$ on $S^{2n-1}$.
Moreover, the centralizer $\mathcal C_{\mathrm{PU}(n,1)}(\mu(\RR^+))$
of $\mu(\RR^+)$ in $\mathrm{PU}(n,1)$ is
one of the following groups up to conjugacy:

\begin{equation}\label{hyperbolic-action}
\mathcal R\times \mbox{\rm U}(n-1)\ \ \mbox{or}\ \
\RR^+\times \mbox{\rm U}(n-1).
\end{equation}
Since $\pi_1(M)$ centralizes $\RR\times \RR^+$,
note that $Q$ centralizes $\RR^+$ (cf. \eqref{central}).
The holonomy group $\mu(Q)$ belongs to
$\mathcal C_{\mathrm{PU}(n,1)}(\mu(\RR^+))$.
As $\mu(Q)$ is a finite subgroup, \eqref{hyperbolic-action} implies that
\begin{equation}\label{finite-center}
\mu(Q)\subset{\rm U}(n-1).
\end{equation}

\subsection{Rigidity of $(M,g,J)$ under the $LCR$ action of $\RR^+$.}
Let $(\eta_0,J_0)$ be the standard 
strictly pseudoconvex 
pseudo-Hermitian structure on $S^{2n-1}$ (cf. \eqref{stand}).
By definition, there exists a positive function $u$ on
$W$ such that
\begin{equation}\label{u-positive}
{\dev}^*\eta_0=u\cdot \eta.
\end{equation}
By Lemma \ref{re}, we know that $A$ is the characteristic
$CR$-vector field on $W$ for $(\eta,J)$.
 If $\{\psi'_t\}$ is
the flow generated by $A$, then note from
\eqref{crtrans} that
$\{{\psi'}_t\}\subset{\rm PSH}(W,\eta,J)$.
Because $W$ is compact, ${\rm PSH}(W,\eta,J)$ is compact.
As ${\rm PSH}(W,\eta,J)\subset \mathop{\Aut}_{CR}(W)$,
the closure of the holonomy image $\mu(\{{\psi'}_t\})$
(which is a connected abelian group) lies in the
maximal torus $T^n$ of the maximal compact subgroup
${\rm U}(n)$ in $\mbox{\rm PU}(n,1)$
up to conjugacy. 
We can describe it as
\[
\mu({\psi'}_t)=(e^{ia_1\cdot t},\cdots,e^{ia_n\cdot t})\ \ \ (\forall t\in\RR)
\]for some $a_i\in\RR$ $(i=1,\ldots,n)$. 
On the other hand, let $\mathcal A=\mathop{\dev}_*(A)$.
Since $\mathop{\dev}$ is equivariant,
$\mathop{\dev}({\psi'}_t w)=\mu({\psi'}_t){\dev}(w)$
on\\
$S^{2n-1}=\{z=(z_1,z_2,\cdots,z_n)\in \CC^n\ |\
|z_1|^2+|z_2|^2+\cdots+|z_n|^2=1\}$, we have:

\begin{equation}\label{vector A}
\mathcal A_z=\frac{ d\mu({\psi'}_t)}{dt}=
\sum_{j=1}^{n}a_j(x_j\frac d{dy_j}\ -\ y_j\frac d{dx_j}) 
\ \ (z=\mathop{\dev}(w), \ \ z_j=x_j+{\mathrm i}y_j).
\end{equation}
As $\eta(A)=1$, we have

\begin{equation}\label{uform}
u(w)={\dev}^*\eta_0(A)=\eta_0(\mathcal A_z)= \sum_{j=1}^{n}a_j\cdot|z_j|^2.
\end{equation}
Since $u>0$ from \eqref{u-positive}, we can assume that
\begin{equation}\label{ai}
0<a_1\leq\cdots\leq a_n.
\end{equation}
As $\mathop{\dev}^{-1}$ maps
the pseudo-Hermitain structure $(\eta,J)$ on $W$
to $(\mathop{\dev}^{-1*}\eta,J_0)$ on $S^{2n-1}$, we put

\begin{equation}\label{pullback form}
\eta_{\mathcal A}={\dev}^{-1*}\eta.
\end{equation} Using \eqref{uform},
we obtain:
\begin{equation}\label{newcon}
\eta_{\mathcal A}=\frac{1}{\mathop{\sum}_{j=1}^{n}a_j\cdot|z_j|^2}\cdot
\eta_0 \ \ \mbox{on}\ S^{2n-1}.
\end{equation} When we note that
$\eta_0=u'\cdot \eta_{\mathcal A}$ where 
$u'=u\circ \mathop{\dev}^{-1}$, and
$\displaystyle T(\RR\times S^{2n-1})=
\{\frac d{dt},\mathcal A\}\oplus\mathop{\Null}\ \eta_0$,
denote the complex structure $J_{\mathcal A}$ 
on $\RR\times S^{2n-1}$ by
\begin{equation}\label{newcomplex}
\begin{split}
&J_{\mathcal A}\frac{d}{dt}=-\mathcal A,\ \ J_{\mathcal A} \mathcal A=
\frac{d}{dt}\\
&J_{\mathcal A}|{\Null}\ \eta_0=J_0.
\end{split}
\end{equation}(Compare \S 3.)
Let ${\rm Pr}:\RR\times S^{2n-1}\ra S^{2n-1}$ be
the canonical projection.
In view of \eqref{lckpara}, setting
\begin{equation}\label{forms}
\begin{split}
&\Omega_{\mathcal A}=d(e^t\cdot {\rm Pr}^*\eta_{\mathcal A}),\ \
\tilde \omega_{\mathcal A}=2e^{-t}\cdot \Omega_{\mathcal A},\\
& \tilde g_{\mathcal A}(X,Y)=\tilde \omega_{\mathcal A}(J_{\mathcal A}X,Y),
\end{split}
\end{equation}we obtain a l.c.K. structure
$(\Omega_{\mathcal A}, J_{\mathcal A})$
on $\RR\times S^{2n-1}$ endowed with the group\\ 
$\RR\times {\rm PSH}(S^{2n-1},\eta_{\mathcal A},J_0)$
of holomorphic homothetic transformations.

\begin{pr}\label{equivariant iso}
There exists an equivariant holomorphic isometry between\\
$(\mathcal C_{\mathcal H}(\RR),\tilde M,\Omega,J)$ and 
$(\RR\times {\rm PSH}(S^{2n-1},\eta_{\mathcal A},J_0), \RR\times S^{2n-1},
\Omega_{\mathcal A},J_{\mathcal A})$.
\end{pr}

\begin{proof}
Let $G:\tilde M\ra \RR\times S^{2n-1}$
be a diffeomorphism defined by
$G(\f_tw)=(t,\mathop{\dev}(w))$.
Note that
${\rm Pr}\circ G=\mathop{\dev}\circ \pi$ on $\tilde M$.
As every element of $\mathcal C_{\mathcal H}(\RR)$
is described as $\f_s\cdot q(\al)$
from Remark \ref{Seifertaction},
define a homomorphism $\Psi:\mathcal C_{\mathcal H}(\RR)\ra
\RR\times {\rm PSH}(S^{2n-1},\eta_{\mathcal A},J_0)$ by setting
\[
\Psi(\f_s\cdot q(\al))=(s,\mu(\al)).
\]
Recall that the action $q(\al)(\f_tw)=\f_t\al w$
from \eqref{split-action}. Then,
\begin{equation*}
\begin{split}
&G(\f_s\cdot q(\al)(\f_tw))
=G(\f_{s+t}\cdot \al w)=(s+t,\mathop{\dev}(\al w))\\
&=(s+t,\mu(\al)\mathop{\dev}(w))=
(s,\mu(\al))(t,\mathop{\dev}(w))=
\Psi(\f_s\cdot q(\al))G(\f_t w).
\end{split}
\end{equation*}
Hence, $(\Psi,G):(\mathcal C_{\mathcal H}(\RR),\tilde M)\ra 
(\RR\times {\rm PSH}(S^{2n-1},\eta_{\mathcal A},J_0),
\RR\times S^{2n-1})$ is equivariantly diffeomorphic.
Next, since $G^*t=t$ for the $t$-coordinate of $\RR\times S^{2n-1}$
and $\mathop{\dev}^*\eta_{\mathcal A}=\eta$ from \eqref{pullback form},
it follows that:
\begin{equation}\label{isometry}
G^*\Omega_{\mathcal A}=
G^*d(e^t\cdot {\rm Pr}^*\eta_{\mathcal A})
=d(e^{G^*t}\cdot G^*{\rm Pr}^*\eta_{\mathcal A})
=d(e^t\cdot \pi^*\eta)=\Omega.
\end{equation}
By definition, $\displaystyle G_*\xi=\frac d{dt}$. Moreover,
when $x=\f_s w$,
\begin{equation*}
G(\psi_t(x))=G(\f_s\psi_t w)=G(\f_s i{\psi'}_t w)=
(s,\mathop{\dev}({\psi'}_t w))=(s,\mu({\psi'}_t)\mathop{\dev}(w)).
\end{equation*}Using \eqref{diff-vec} and 
\eqref{vector A},
\begin{equation*}
G_*(-J\xi_x)=\frac{dG\psi_t}{dt}(x)|_{t=0}=\mathcal A_{\small{Gx}}
=-J_{\mathcal A}(\frac d{dt})_{\small{Gx}}.
\end{equation*}
Thus $G_*(J\xi)=J_{\mathcal A}G_*\xi$.
As $G^*\Omega_{\mathcal A}=\Omega$ from \eqref{isometry},
$G$ maps
$\{\xi,J\xi\}^{\perp}$ onto
$\displaystyle \{\frac d{dt},\mathcal A\}^{\perp}$.
Consider the commutative diagram:
\begin{equation}\label{JA-com}
\begin{CD}
(\{\xi,J\xi\}^\perp,J)@>\pi_*>>({\Null}\, \eta,J)\\
@VV{G_*} V                 @VV{\dev_*} V\\
(\{\frac d{dt},\mathcal A\}^\perp,J_{\mathcal A})@>{\rm Pr}_*>>({\Null}\,
 \eta_0,J_0).
\end{CD}
\end{equation}Here note that $J_{\mathcal A}=J_0$
on ${\Null}\, \eta_{\mathcal A}={\Null}\, \eta_0$.
For $X\in\{\xi,J\xi\}^\perp$,
\begin{equation*}
{\rm Pr}_*G_*J(X)=
{\dev}_*(J\pi_*X)=J_0{\dev}_*\pi_*(X)=J_{\mathcal A}{\rm Pr}_*G_*(X)=
{\rm Pr}_*J_{\mathcal A}G_*(X),
\end{equation*}thus,
$G_*J(X)=J_{\mathcal A}G_*(X)$.
Hence, $G$ is $(J,J_{\mathcal A})$-biholomorphic.
Moreover, as\\
$G^*\tilde\omega_{\mathcal A}
=G^*(2e^{-t}\Omega_{\mathcal A})=2e^{-t}\Omega=\bar \Theta$
and $\bar g(X,Y)=\bar \Theta(JX,Y)$,
we obtain that $G^*\tilde g_{\mathcal A}=\bar g$.
Therefore, $(\Psi,G)$ induces a holomorphic isometry
from $(M,\hat g,J)$ onto\\
$(\RR\times S^{2n-1}/\Psi(\pi_1(M)),
\hat g_{\mathcal A},\hat J_{\mathcal A})$.

\end{proof}

\subsection{The Hopf manifold 
${\bf \RR\times S^{2n-1}/\Psi(\pi_1(M))}$}
We prove that
$\RR\times S^{2n-1}/\Psi(\pi_1(M))$
is a primary Hopf manifold $M_\Lambda$ for some $\Lambda$
obtained in \S 3.
Each element of $\pi_1(M)$
is of the form $\gamma=\f_s\cdot q(\al)$ for some $s\in\RR$
where $\nu(\gamma)=\al\in Q=\nu(\pi_1(M))$.
By definition of $\Psi$,
$\Psi(\gamma)=(s,\mu(\al))$.
We show that $\Psi(\pi_1(M))$ has no torsion element.
For this, if $\Psi(\gamma)$
is of finite order (say, $\ell$),
then $1=(0,1)=\Psi(\gamma^{\ell})=(\ell s,\mu(\al^{\ell}))$.
Then, $s=0$ so that $\Psi(\gamma)=(0,\mu(\al))$.
On the other hand, recall from \eqref{finite-center} that
$\mu(Q)\subset {\rm U}(n-1)$ up to conjugacy, and so
$\mu(Q)$ has a fixed point $w_0\in S^{2n-1}$.
Since $\Psi(\pi_1(M))$ acts freely on $\RR\times S^{2n-1}$, 
while $\Psi(\gamma)(t,w_0)=(t,\mu(\al)w_0)=(t,w_0)$,
it follows that $\Psi(\gamma)=1$.
Moreover, if $\gamma_1=\f_{s_1}\cdot q(\al_1)$,
$\gamma_2=\f_{s_2}\cdot q(\al_2)$, then
$\Psi([\ga_1,\ga_2])=(0,\mu([\al_1,\al_2])$. By the same reason,
$\Psi([\pi_1(M),\pi_1(M)])=\{1\}$.
Hence, $\pi_1(M)$ is a finitely generated torsionfree abelian group. 
If we recall from \eqref{central}
that $1\ra \ZZ\ra \pi_1(M)\stackrel{\nu}{\lra} Q \ra 1$
is the central group extension where $Q$ is finite, then
$\pi_1(M)$ itself is an infinite cyclic group. 
Since
$\Psi(\pi_1(M))\subset\RR\times {\rm PSH}(S^{2n-1},\eta_{\mathcal A},J_0)$
and the projection maps $\Psi(\pi_1(M))$ onto 
$\mu(Q)$ in ${\rm PSH}(S^{2n-1},\eta_{\mathcal A},J_0)$,
$\mu(Q)$ is a finite cyclic group.
As ${\rm PSH}(S^{2n-1},\eta_{\mathcal A},J_0)$ has the maximal torus $T^n$
(cf. \eqref{maximaltorus}),
we obtain that $\Psi(\pi_1(M))\subset \RR\times T^n$
up to conjugacy.
A generator of $\Psi(\pi_1(M))$ is described as
 $(s,(c_1,\cdots,c_n))\in \RR\times T^n$. 
Noting \eqref{ai}, let
$\lambda_j=e^{-a_js}c_j$ and $\Lambda=(\lambda_1,\cdots,\lambda_n)$.
By Theorem \ref{secondary} and the remark below,
$\RR\times S^{2n-1}/\Psi(\pi_1(M))$ is a primary Hopf manifold
 $M_{\Lambda}$ of type $\Lambda$.
This finishes the proof of Theorem C in the Introduction.\\[1mm]


\begin{thebibliography}{99}
\bibitem{Be} F.A. Belgun, \emph{On the metric structure of non-K\"ahler complex
 surfaces}, Math. Ann., {\bf 317} (2000), 1-40.
\bibitem{Bl} D.E. Blair, Contact manifolds in Riemannian geometry,
L.N.M. {\bf 509}, Springer Verlag 1976.
\bibitem{CG} S.S. Chern and L. Greenberg, \emph{Hyperbolic Spaces},
in \emph{Contribution to Analysis} (A Collection of papers Dedicated to 
Lipman Bers, eds. L. Ahlfors and others), Academic Press, New York and London,
49-87, 1974.
\bibitem{CR} P. Conner and F. Raymond,
\emph{Injective operation of the toral groups}, Topology, {\bf 10} (1971),
283-296.
\bibitem {DO} S. Dragomir and L. Ornea, Locally conformal
K{\"a}hler geometry, Progress in Math. {\bf 155}, Birkh{\"a}user 1998.
\bibitem {GO} P. Gauduchon and L. Ornea, \emph{Locally conformally K\"ahler
metrics on Hopf surfaces},
Annales de l'Inst. Fourier, {\bf 48} (1998), 1107-1127.
\bibitem{Gold} W. Goldman, \emph{Complex hyperbolic geometry},
 Oxford Mathematical Monographs,
Oxford Univ. Press, 1999.
\bibitem{Kam} Y. Kamishima, \emph{Geometric flows on compact manifolds
and global rigidity}, Topology, {\bf 35} (1996), 439-450.
\bibitem{kamicomp} Y. Kamishima, \emph{Holomorphic torus actions on compact
locally conformal K\"ahler manifolds},
Compositio Math., {\bf 124} (2000), 341-349.
\bibitem {Ko} S. Kobayashi, Transformation groups in differential
geometry, Ergebnisse der Math. {\bf 70}, Springer Verlag, 1972.
\bibitem {Ko1} S. Kobayashi and K. Nomizu,
Foundationsv of differential geometry II, Interscience Publisheres, New York,
1969.
\bibitem{Ku} R. Kulkarni, \emph{On the principle of uniformization},
J. Diff. Geom., {\bf 13} (1978), 109-138.
\bibitem{Le} J. Lelong-Ferrand,
\emph{Transformations conformes et quasi conformes des
vari\'et\'es  riemanniennes compactes}, Acad. Roy. Belgique Sci. Mem. Coll., 
{\bf 8} (1971), 1-44.
\bibitem{Ob} M. Obata,
\emph{The conjectures on conformal transformations of Riemannian manifolds},
J. Diff. Geom., {\bf 6} (1971), 247-258.
\bibitem{Ra} M. S. Raghunathan, Discrete subgroups of Lie groups,
Springer Verlag, 1972.
\bibitem{Tr} F. Tricerri, \emph{Some examples of locally conformal
K\"ahler manifolds}, Rend. Sem. Mat. Univ. Politecn. Torino {\bf 40} (1982),
81-92.
\bibitem {Va} I. Vaisman, {\em Locally conformal
K\"ahler manifolds with parallel Lee form}, Rend. Mat. {\bf 12} (1979),
263-284.
\bibitem {Va1} I. Vaisman, {\em Generalized Hopf manifolds},
Geometriae Dedicata, 13(1982), 231-255.
\bibitem{W} S. M. Webster,
\emph{On the transformation group of a real hypersurface}, Trans. Amer. Math.
Soc., {\bf 231} (1977), 179-190.
\bibitem{W1} S. M. Webster,
\emph{Pseudohermitian geometry of a real hypersurface}, J. Diff. Geom., 
{\bf 13} (1978), 25-41.

\end{thebibliography}
\end{document}